\documentclass[11pt,letterpaper]{amsart}
\usepackage{amsmath,amssymb}
\usepackage[T1]{fontenc}
\usepackage{lmodern}
\usepackage{hyperref}
\hypersetup{hidelinks,pdfauthor={Hongrui Yuan},pdftitle={The Casimir on the K-fixed tensor square of the minimal representation of G2(2)}}
\newcommand{\R}{\mathbb R}
\newcommand{\C}{\mathbb C}
\newcommand{\N}{\mathbb N}
\newcommand{\g}{\mathfrak g}
\newcommand{\kk}{\mathfrak k}
\newcommand{\pp}{\mathfrak p}
\newcommand{\Cas}{\operatorname{Cas}}
\newcommand{\Dom}{\operatorname{Dom}}
\newcommand{\spfin}{\operatorname{span}_{\mathrm{fin}}}
\newcommand{\dd}{\,\mathrm d}
\newcommand{\rbra}[1]{\left(#1\right)}
\newcommand{\fs}[2]{\frac{#1}{#2}}
\newcommand{\coloneqq}{\mathrel{:=}}
\newcommand{\U}{\mathcal U}

\newtheorem{theorem}{Theorem}[section]

\newtheorem{lemma}[theorem]{Lemma}
\newtheorem{corollary}[theorem]{Corollary}
\newtheorem{conjecture}[theorem]{Conjecture}
\theoremstyle{remark}
\newtheorem{remark}[theorem]{Remark}
\numberwithin{equation}{section}
\title[Casimir on the $K$-fixed tensor square]{The Casimir on the $K$-fixed tensor square of the minimal representation of $G_{2(2)}$}
\author{Hongrui Yuan}
\email{ytyh@pku.edu.cn}
\date{Oct 16, 2026}
\begin{document}
\begin{abstract}
We compute the quadratic Casimir on the $K$-fixed subspace of the Hilbert tensor square of the unitary minimal representation of the split real group $G_{2(2)}$. In an explicit algebraic basis of invariant tensors, the Casimir acts by a three-term recurrence. Normalizing this basis yields a symmetric Jacobi operator whose self-adjoint realization is unitarily equivalent to multiplication by $-\frac92-8x^2$ in a Wilson-polynomial model with parameters $(\frac34,\frac34,\frac7{12},\frac{11}{12})$. We give the normalized spectral density and the unitary transform. The resulting spectrum is purely absolutely continuous, equals $(-\infty,-\frac92]$, and has scalar spectral multiplicity one.
\end{abstract}
\maketitle
\section{Introduction}
The decomposition of a tensor product of unitary representations is a basic problem in representation theory. Even when both factors are explicitly realized, the diagonal action can be difficult to describe and analyze. Central elements of the enveloping algebra provide accessible spectral data: their action can often be computed on a subspace fixed by a maximal compact subgroup, where compact representation theory reduces the problem to a sequence of finite-dimensional calculations.

We consider the connected linear split real group $G=G_{2(2)}$ and its unitary minimal representation $\pi_{\min}$ on $\mathcal H_{\min}$. Existence and unitarity are known; Gross--Wallach \cite[Proposition 14.2]{wallach_quaternionic_1996} record Vogan's construction. The conformally invariant systems of Barchini--Kable--Zierau \cite{barchini_conformally_2008,barchini_conformally_2009} and Kable's description of their $K$-finite solutions \cite{kable_k-finite_2011,kable_conformally_2012} supply the geometric model, compact-type ladder, embedding vectors, and noncompact action used here. These inputs make it possible to carry out the tensor calculation in fixed compact coordinates.

The representation and its one-factor action are established inputs. The normalization of the invariant form and the construction of compact invariant tensors support the tensor calculation. Frahm, van Haastrecht, Weiske, and Zhang \cite[Theorem 3.6 and Proposition 4.5]{frahm_degenerate_2026} describe the same module and form through a degenerate induced family and intertwiners. Here we compute the tensor-square Casimir in the compact model and obtain its representation-specific Wilson realization from the established polynomial theory \cite[Section 9.1]{koekoek_hypergeometric_2010}. 

Let $K$ be a maximal compact subgroup of $G$. The precise object studied in this article is the quadratic Casimir for the diagonal representation $\Pi=\pi_{\min}\otimes\pi_{\min}$ on
\[
\mathcal H^K=(\mathcal H_{\min}\widehat\otimes\mathcal H_{\min})^K.
\]
The compact-type decomposition is multiplicity free, and compact self-duality shows that only matching types contribute to $\mathcal H^K$. Each such pair contains the one-dimensional space of $K$-invariant vectors. Choosing its algebraic generator $\Psi_n$ explicitly allows us to compute the Casimir before making any spectral normalization.

Our main algebraic result is the recurrence in Theorem~\ref{thm:4.cas_action}. For the invariant bilinear form $B=\frac18B_{\mathfrak g}$, where $B_{\mathfrak g}$ is the Killing form, it reads
\[
\mathrm d\Pi(\Cas)\Psi_n
=-2\left(n+\frac53\right)^2\Psi_{n+1}
-2\delta_n\Psi_n
-2\left(n+\frac13\right)^2\Psi_{n-1},
\]
where $\Psi_{-1}=0$ and
\[
\delta_n=\frac{(n+2)(n+\frac43)(n+\frac53)+n(n+\frac13)(n+\frac23)}{n+1}.
\]
The adjacent one-factor action explains the tridiagonal structure. Applying the noncompact elements successively to the explicit invariant tensors determines all three coefficients and verifies the cancellations between output compact types.

Our main spectral result, Theorem~\ref{thm:4.spectral_decomp}, identifies the self-adjoint Casimir on $\mathcal H^K$ with a scalar multiplication operator. After normalizing $\Psi_n$ to an orthonormal basis $\Phi_n$, we obtain the Wilson parameters
\[
(a,b,c,d)=\left(\frac34,\frac34,\frac7{12},\frac{11}{12}\right).
\]
The unitary transform sends $\Phi_n$ to $(n+1)(-1)^n\widetilde W_n(x^2;a,b,c,d)$ and the Casimir to multiplication by $-\frac92-8x^2$. The normalized spectral density is explicit. Consequently the spectrum is purely absolutely continuous, equals $(-\infty,-\frac92]$, has scalar spectral multiplicity one, and has no discrete atoms.

The layout of this document is organized as follows. Section~\ref{sec:preliminaries} fixes the special-function and compact conventions. Section~\ref{sec:model} describes the geometric model, and Section~\ref{sec:kfinite} passes from its $K$-finite vectors to the noncompact action and invariant form. Section~\ref{sec:basis} constructs the invariant tensors and computes the Casimir. Section~\ref{sec:spectral} normalizes the basis and proves the spectral realization. Section~\ref{sec:conjectural} records a conjectural spherical decomposition suggested by this spectral model. The appendices give, in that order, the compact normalization and action calculations, followed by the tensor contractions.

\section{Preliminaries}\label{sec:preliminaries}
We collect the special-function formulas and compact-representation conventions that will be used in the calculations. All Hermitian inner products are linear in the first variable.

\subsection{Wilson polynomials}
For a nonnegative integer $j$, the rising factorial of the Pochhammer symbol is defined by $(z)_j=z(z+1)\cdots(z+j-1)$, with $(z)_0=1$. We abbreviate products by $(a_1,\ldots,a_r)_j=\prod_{i=1}^r(a_i)_j$. When an upper parameter is a nonpositive integer $-n$, the hypergeometric notation denotes the finite sum
\[
{}_rF_s\!\left(\begin{matrix}-n,a_2,\ldots,a_r\\b_1,\ldots,b_s\end{matrix};z\right)
=\sum_{j=0}^{n}\frac{(-n,a_2,\ldots,a_r)_j}{(b_1,\ldots,b_s)_j}\frac{z^j}{j!},
\]
provided the denominator factors do not vanish.

For $a,b,c,d>0$, write $S=a+b+c+d$ and use the normalization
\begin{equation}\label{eq:Wilson-definition}
\begin{aligned}
\widetilde W_n(x^2;a,b,c,d)
&=\frac{W_n(x^2;a,b,c,d)}{(a+b,a+c,a+d)_n}\\
&={}_4F_3\!\left(\begin{matrix}-n,n+S-1,a+ix,a-ix\\a+b,a+c,a+d\end{matrix};1\right).
\end{aligned}
\end{equation}

The Wilson polynomials are orthogonal for the positive measure
\begin{equation}\label{eq:Wilson-weight}
\omega_W(x)\dd x
=\frac{1}{2\pi}\left|\frac{\Gamma(a+ix)\Gamma(b+ix)\Gamma(c+ix)\Gamma(d+ix)}{\Gamma(2ix)}\right|^2\dd x,
\quad x>0,
\end{equation}
where $\Gamma(z)$ is the Euler gamma function; the orthogonality relation and the corresponding squared norms for this normalization are given in \cite[Section 9.1]{koekoek_hypergeometric_2010}.
\begin{lemma}[Wilson completeness]\label{lem:Wilson-complete}
For $a,b,c,d>0$, the normalized Wilson polynomials form a complete orthogonal system in $L^2((0,\infty),\omega_W(x)\dd x)$.
\end{lemma}
\begin{proof}
The positive-parameter Wilson measure is determinate, so the polynomials are dense in the full $L^2$ space. This is the completeness result in \cite[Section 9.1]{koornwinder2025additionsformulalistshypergeometric}, with the determinacy and density criteria of \cite[Theorems 2.9 and 2.14]{shohat1943problem}. Multiplying each polynomial by its nonzero normalization factor does not change its span.
\end{proof}

For $n\geq1$, the three-term recurrence for the normalized Wilson polynomials is
\begin{equation}\label{eq:Wilson-recurrence}
-(a^2+x^2)\widetilde W_n(x^2)
=A_n\widetilde W_{n+1}(x^2)-(A_n+C_n)\widetilde W_n(x^2)+C_n\widetilde W_{n-1}(x^2),
\end{equation}
where, for $n\geq1$,
\begin{equation}\label{eq:Wilson-coefficients}
\begin{aligned}
A_n&=\frac{(n+S-1)(n+a+b)(n+a+c)(n+a+d)}{(2n+S-1)(2n+S)},\\
C_n&=\frac{n(n+b+c-1)(n+b+d-1)(n+c+d-1)}{(2n+S-2)(2n+S-1)}.
\end{aligned}
\end{equation}
These formulas are from \cite[Section 9.1]{koekoek_hypergeometric_2010}. At $n=0$, the terminating definition \eqref{eq:Wilson-definition} gives separately
\[
\widetilde W_0=1,\quad
\widetilde W_1=1-\frac{S(a+ix)(a-ix)}{(a+b)(a+c)(a+d)}
=1-\frac{S(a^2+x^2)}{(a+b)(a+c)(a+d)}.
\]
Consequently, with $A_0=\frac{(a+b)(a+c)(a+d)}{S}$, the initial relation is
\begin{equation}\label{eq:Wilson-initial}
-(a^2+x^2)\widetilde W_0=A_0(\widetilde W_1-\widetilde W_0).
\end{equation}
This uses only the constant polynomial and the two terms of the degree-one polynomial; no polynomial with index $-1$ is introduced. We set $C_0=0$ when recording the recurrence coefficients.

\subsection{Compact modules and Clebsch--Gordan maps}
For compact calculations, let $\Gamma_q$ be the irreducible $\mathrm{SU}(2)$-module of highest weight $q$ and dimension $q+1$, where $q$ is a nonnegative integer.

Its weight basis is $\{\xi^q_{q-2j}:0\leq j\leq q\}$; the subscript specifies the weight. 

Let $A_1,A_2,A_3$ be the standard real compact $\mathfrak{su}(2)$-basis, with $[A_1,A_2]=2A_3$, $[A_2,A_3]=2A_1$, and $[A_3,A_1]=2A_2$. For convenience, set $A_\pm=\pm A_2-iA_3$. Their useful action formulas are
\begin{equation}\label{eq:compact-action}
A_1\xi_a^q=ia\xi_a^q,\quad
A_+\xi_a^q=(q-a)\xi_{a+2}^q,\quad
A_-\xi_a^q=(q+a)\xi_{a-2}^q.
\end{equation}
The invariant Hermitian form is fixed by
\begin{equation}\label{eq:compact-norm}
\langle\xi^q_{q-2j},\xi^q_{q-2k}\rangle_q
=\binom qj^{-1}\delta_{jk}.
\end{equation}
Thus a highest vector has norm one. These compact forms specify the meaning of an isometric tensor-product decomposition.

The Clebsch--Gordan decomposition is
\[
\Gamma_n\otimes\Gamma_m\simeq
\Gamma_{n+m}\oplus\Gamma_{n+m-2}\oplus\cdots\oplus\Gamma_{n-m}
\quad(n\geq m),
\]
with multiplicity one; see \cite[Section 8.1.3]{vilenkin_representation_2012}.
We use in particular the unitary maps
\begin{equation}\label{eq:CG-interface}
\begin{aligned}
&f_m\colon\Gamma_m\otimes\Gamma_1\longrightarrow\Gamma_{m+1}\oplus\Gamma_{m-1},\\
&g_n\colon\Gamma_n\otimes\Gamma_3\longrightarrow
\Gamma_{n+3}\oplus\Gamma_{n+1}\oplus\Gamma_{n-1}\oplus\Gamma_{n-3},
\end{aligned}
\end{equation}
in the metrics \eqref{eq:compact-norm}. Their phases are fixed by requiring every nonzero component of $f_m(\xi_m^m\otimes\xi_{-1}^1)$ and $g_n(\xi_n^n\otimes\xi_{-3}^3)$ to have positive coefficient in the displayed weight bases. In particular, the upper projection sends the tensor of highest vectors to the highest vector with coefficient one. These are the normalized maps obtained from the coefficients in \cite[Section 8.1.3]{vilenkin_representation_2012} by conversion from orthonormal bases to \eqref{eq:compact-norm}. The projection formulas and normalization checks used below are recorded in Appendix~\ref{app:compact}.

\begin{remark}
In low dimensions only the summands in the actual compact product are used: $\Gamma_0\otimes\Gamma_1=\Gamma_1$ and $\Gamma_2\otimes\Gamma_3=\Gamma_5\oplus\Gamma_3\oplus\Gamma_1$. A nonexistent summand is omitted before any coefficient is evaluated. Within an existing summand, a vector outside its weight range is zero. These rules keep the lowest-level calculations separate from generic formulas with vanishing denominators.
\end{remark}

\subsection{The split group and its inducing data}
The group conventions specify both the compact factors and the parabolic used for the geometric model. We follow \cite[Section 2]{kable_conformally_2012}. Let $\g$ be the complex Lie algebra of type $G_2$, with short simple root $\alpha$, long simple root $\beta$, and
\[
(\alpha,\alpha)=\frac23,\quad(\beta,\beta)=2,\quad(\alpha,\beta)=-1.
\]
Fix the Chevalley basis of that reference, with $H_\mu=[X_\mu,X_{-\mu}]$, and use the invariant bilinear form
\begin{equation}\label{eq:B}
B=\frac18B_{\g}\iff B(X_\mu,X_{-\mu})=\frac{2}{(\mu,\mu)},
\end{equation}
where $B_{\g}$ is the Killing form. Write $\vartheta=3\alpha+2\beta=\varpi_\beta$ for the highest root. The eigenspaces of $\operatorname{ad}H_\vartheta$ give the grading $\g=\bigoplus_{j=-2}^{2}\g(j)$ and the parabolic subalgebra $\mathfrak q=\bigoplus_{j=0}^{2}\g(j)$.

The real span $\g_0$ of the Chevalley basis is the split real form. Let $G=G_{2(2)}$ be the connected linear split real group with this Lie algebra. The Cartan involution sends $X_\mu$ to $-X_{-\mu}$ and $H$ to $-H$, giving $\g_0=\kk_0\oplus\pp_0$. Its maximal compact subgroup is
\[
K\simeq(\mathrm{SU}(2)\times\mathrm{SU}(2))/\{(1,1),(-1,-1)\}.
\]
The two commuting compact triples, in the order denoted $U$ and $V$ in \cite[Section 2]{kable_conformally_2012}, identify $\kk_0$ with $\mathfrak{su}(2)\oplus\mathfrak{su}(2)$. For modules of these two factors, $\boxtimes$ denotes the external tensor product: the first factor acts on the first module and the second on the second. We keep this compact-factor order throughout.

Let $\pp$ be the complexification of $\pp_0$. By \cite[Lemma 2.2]{kable_conformally_2012}, its compact-type identification is isometric for the Hermitian form $\langle Y,Z\rangle_{\pp}=B(Y,\bar Z)$ and the product of the forms \eqref{eq:compact-norm}. We write it in the direction
\begin{equation}\label{eq:iota}
\iota\colon\pp\longrightarrow\Gamma_3\boxtimes\Gamma_1,
\quad
\iota\left(\frac{H_\vartheta+i(X_\vartheta+X_{-\vartheta})}{2}\right)
=\xi_3^3\boxtimes\xi_1^1.
\end{equation}
Our $\iota$ is the inverse of the map displayed there, with the same scale and phases. The bar denotes conjugation relative to $\g_0$; in particular $\overline{X+iZ}=X-iZ$ for real $X,Z$.

Let $Q=LN$ be the real Heisenberg parabolic with complexified Lie algebra $\mathfrak q$, where $L$ has complexified Lie algebra $\g(0)$ and $N$ is the unipotent radical; write $\bar Q=L\bar N$ for the opposite parabolic. The Levi factor is isomorphic to $\mathrm{GL}(2,\R)$. Its real analytic characters are parameterized by $s\in\R$ and $\varepsilon\in\{0,1\}$ through
\begin{equation}\label{eq:characters}
\mathrm d\chi(s,\varepsilon)=s\varpi_\beta,
\quad
\chi(s,\varepsilon)(\operatorname{diag}(1,-1))=(-1)^\varepsilon.
\end{equation}
The continuous parameter fixes the differential of the inducing character, while $\varepsilon$ fixes its value on the other component of $L$. These are the parameters used in the model below, in Kable's convention.

\section{The geometric realization}\label{sec:model}
The tensor-product calculation starts from a realization of one minimal representation in which compact vectors and infinitesimal actions can be computed. We use the known geometric model from \cite{barchini_conformally_2008,barchini_conformally_2009,kable_conformally_2012}. The purpose of this section is to specify the ambient space and the equations selecting that representation.

\subsection{The ambient induced representation}
For $\eta=\chi(s,\varepsilon)$, extend $\eta$ to $\bar Q$ trivially on $\bar N$. The homogeneous line bundle $\mathcal L_\eta$ on $G/\bar Q$ is defined by
\[
(g\bar q,z)\sim(g,\eta(\bar q)^{-1}z).
\]
Its smooth sections are represented by functions $F:G\to\C$ satisfying
\begin{equation}\label{eq:covariance}
F(g\bar q)=\eta(\bar q)F(g),
\quad g\in G,\quad\bar q\in\bar Q.
\end{equation}
Left translation, $\pi_\eta(g)F(x)=F(g^{-1}x)$, makes this the unnormalized induced representation $\operatorname{Ind}_{\bar Q}^{G}\eta^{-1}$. The ambient induced space contains more than the representation we seek; a differential system will select the relevant subspace.

\subsection{The differential system and the parameter choice}
To describe the equations, restrict first to the open cell $U_0=N\bar Q/\bar Q$, identified with $N$. In this chart, the algebraic construction of differential operators can be described by right differentiation. If $\U$ denotes the universal enveloping algebra, elements of $\U(\mathfrak n)$ give operators commuting with the left $N$-action, with
\[
(D_XF)(n)=\left.\frac{\mathrm d}{\mathrm dt}F(n\exp(tX))\right|_{t=0}
\quad(X\in\mathfrak n_0).
\]
The formulation in the generalized Verma module incorporates the inducing character and gives the corresponding operators on the line bundle; see \cite[pp.\,800--801]{barchini_conformally_2009} and \cite[Section 3]{kable_conformally_2012}. We need this correspondence to specify the equations and, subsequently, to turn them into conditions on compact vectors.

The equations must be compatible with the infinitesimal $G$-action if their common solution space is to carry a representation. This is the role of conformal invariance. A system $D_1,\ldots,D_r$ is conformally invariant when, for each $Y\in\g$,
\begin{equation}\label{eq:conformal}
[\mathrm d\pi_\eta(Y),D_i]=\sum_{j=1}^{r}c_{ji}^{Y}D_j
\end{equation}
for smooth coefficient functions $c_{ji}^{Y}$ on the open cell. This relation makes the common kernel stable under the infinitesimal action.

For the Heisenberg parabolic, the relevant system is the second-order system
\[
\Omega_2=\{\Omega_2(X_\alpha),\Omega_2(H_\alpha),\Omega_2(X_{-\alpha})\}.
\]
Its three operators are obtained from the algebraic construction just described, indexed by the three generators of the semisimple part of the Levi algebra. Their explicit formulas and the invariance calculation are given in \cite[Theorem 5.2 and Section 8.9]{barchini_conformally_2008} and \cite[Section 3.2]{kable_conformally_2012}. We use this specified system rather than the general family of systems on a flag manifold.

Conformal invariance holds at $s=\frac23$. There are still two choices of component character, and they do not give the same solution space. Kable proves that the $\Omega_2$-solution space for $\chi(\frac23,0)$ is zero, whereas the one for $\chi(\frac23,1)$ is nonzero and irreducible \cite[Theorem 4.8]{kable_conformally_2012}. We therefore fix
\[
\eta_*:=\chi\!\left(\frac23,1\right),\quad
\mathcal S=\{F\in\Gamma(\mathcal L_{\eta_*}):DF=0\text{ for all }D\in\Omega_2\}.
\]
The symbol $\Gamma(\mathcal L_{\eta_*})$ here denotes global sections, rather than a compact module. This is the geometric realization of the known minimal representation. We write $V_{\min}$ for its $K$-finite module and $\pi_{\min}$ for the action. To make this model explicit, the next task is to find the compact vectors satisfying the equations.

\section{\texorpdfstring{$K$}{K}-finite vectors, noncompact action, and the invariant form}\label{sec:kfinite}

\subsection{Embedding vectors and matrix coefficients}
In this section, we follow Kable's realization \cite{kable_conformally_2012}. The geometric model specifies the representation by equations on sections. For $K$-finite solutions, these equations can be tested in finite-dimensional compact modules. Let $(\sigma,E_\sigma)$ be an irreducible representation of $K$.

An embedding vector must first have the required $K\cap L$ character:
\[
E_\sigma^{(K\cap L,\eta_*)}
=\{u\in E_\sigma:\sigma(l)u=\eta_*(l)u\text{ for all }l\in K\cap L\}.
\]
The differential equations impose an additional condition. Under the compact form of the operator correspondence, let $v_1,v_2,v_3\in\U(\kk)$ represent the three operators of $\Omega_2$. The embedding-vector space is
\begin{equation}\label{eq:embedding-space}
\mathbb M(\sigma)=\{u\in E_\sigma^{(K\cap L,\eta_*)}:
\mathrm d\sigma(\bar v_i)u=0,\ 1\leq i\leq3\}.
\end{equation}
Conjugation on $\U(\kk)$ is induced by the real form $\kk_0$. Thus an embedding vector records both the bundle covariance and the solution equations; it is not an arbitrary vector of $E_\sigma$.

To recover an actual section from such a vector, equip $E_\sigma$ with a $K$-invariant positive Hermitian form and take the matrix coefficient
\begin{equation}\label{eq:matrix-coefficient}
\psi_\sigma(\xi,u)(k)=\langle\xi,\sigma(k)u\rangle_\sigma,
\quad k\in K.
\end{equation}
The variable $\xi$ runs over $E_\sigma$, while $u$ is fixed. The identification $G/\bar Q\simeq K/(K\cap L)$ turns the covariance of this function into the line-bundle condition \eqref{eq:covariance}. Kable's correspondence says that the copies of $\sigma$ in the $K$-finite solution space are precisely those obtained from $u\in\mathbb M(\sigma)$ \cite[Theorem 2.6]{kable_k-finite_2011}. Their multiplicity is $\dim\mathbb M(\sigma)$. This supplies the link from the geometric equations to a compact-type decomposition.

For the selected system, the nonzero embedding spaces occur exactly at
\begin{equation}\label{eq:ladder}
E_m=\Gamma_{3m+2}\boxtimes\Gamma_m,\quad m\geq0,
\quad V_{\min}=\bigoplus_{m\geq0}V_m,\quad V_m\simeq E_m,
\end{equation}
and each embedding space has dimension one \cite[Theorem 4.8]{kable_conformally_2012}. The total highest weight $4m+2$ is even, so these modules descend to $K$. Put $\zeta_a=\xi^{3m+2}_{3a}\boxtimes\xi^m_a$ for $a=m-2k$, $0\leq k\leq m$. The fixed choice of embedding vector is
\begin{equation}\label{equ:4.2}
u_m=\sum_{k=0}^{m}(-1)^k\binom mk
\frac{(\frac43)_m}{(\frac43)_k(\frac43)_{m-k}}\,\zeta_{m-2k}.
\end{equation}
This is Kable's explicit vector, with its leading coefficient one. We use the concrete identification
\[
\varphi_m(\xi)=\psi_{E_m}(\xi,u_m)\in V_m.
\]
In this notation $\psi_{E_m}$ means the matrix coefficient for the compact action on $E_m$, as in \eqref{eq:matrix-coefficient}. The compact structure and the scale of its realized vectors are now fixed.

\subsection{Transferring the noncompact action}
The compact action on \eqref{eq:matrix-coefficient} is already known: it is the action on the first vector $\xi$. To determine the full infinitesimal representation, it remains to compute the noncompact directions $Y\in\pp$. Kable's map $R(\sigma)$ transfers this problem to a finite-dimensional tensor product.

For an orthonormal real basis $\{W_i\}$ of $\pp_0$ with respect to $B$, write $W_i=Z_i+U_i$ with $Z_i\in\kk_0$ and $U_i\in\bar{\mathfrak q}_0=\operatorname{Lie}(\bar Q)$, as in \cite[Section 2]{kable_k-finite_2011}. The transfer map is
\[
R(\sigma)u=\sum_iW_i\otimes
\bigl(\mathrm d\sigma(Z_i)+\mathrm d\eta_*(U_i)\bigr)u,
\quad u\in E_\sigma^{(K\cap L,\eta_*)}.
\]
With the product compact metric on $\pp\otimes E_\sigma$, its relation to infinitesimal action is
\begin{equation}\label{eq:transfer}
\mathrm d\pi_{\min}(Y)\psi_\sigma(\xi,u)
=-\psi_{\pp\otimes\sigma}(Y\otimes\xi,R(\sigma)u),
\end{equation}
for embedding vectors $u$, by \cite[Lemma 2.7]{kable_k-finite_2011}. The same lemma maps embedding vectors to embedding vectors for $\pp\otimes\sigma$. The minus sign in \eqref{eq:transfer} is part of the convention and will also enter the form comparison.

It remains to evaluate $R(\sigma)$ on the vector \eqref{equ:4.2}. Let $R_m=R(E_m)$ and $N=3m+2$. Write $Y_{a,b}=\iota^{-1}(\xi_a^3\boxtimes\xi_b^1)$. The input formula from \cite[Lemma 4.7]{kable_conformally_2012}, at the selected parameter $s=\frac23$, is
\begin{equation}\label{eq:R-zeta}
\begin{aligned}
R_m\zeta_a&=\left(a+\frac23\right)Y_{3,1}\otimes\zeta_a
+\left(a-\frac23\right)Y_{-3,-1}\otimes\zeta_a\\
&\quad+\frac{N-3a}{2}Y_{1,1}\otimes(\xi^N_{3a+2}\boxtimes\xi^m_a)-\frac{N+3a}{2}Y_{-1,-1}\otimes(\xi^N_{3a-2}\boxtimes\xi^m_a)\\
&\quad+\frac{m-a}{2}Y_{3,-1}\otimes(\xi^N_{3a}\boxtimes\xi^m_{a+2})-\frac{m+a}{2}Y_{-3,1}\otimes(\xi^N_{3a}\boxtimes\xi^m_{a-2}).
\end{aligned}
\end{equation}
Thus $R_mu_m$ is obtained by summing this formula with the coefficients of $u_m$.

The result is still expressed in $\pp\otimes E_m$. To identify its compact components, apply $g_n$ and $f_m$ from \eqref{eq:CG-interface}, using the ordinary unsigned tensor flips to place the two compact factors together. The candidate types are
\[
\Gamma_{3m+5-2r}\boxtimes\Gamma_{m+1-2t},
\quad 0\leq r\leq3,\quad0\leq t\leq1,
\]
with absent summands omitted. Individual terms of \eqref{eq:R-zeta} can have components in all these types. Their weighted sum is more restricted: $R$ preserves the embedding-vector condition, and the only candidate types with nonzero embedding spaces are $E_{m+1}$ and $E_{m-1}$. This explains the cancellation of the other six projections for $m\geq1$. The support and the two retained summed coefficients are given in \cite[Theorem 4.8 and its proof, p.\,122]{kable_conformally_2012}.

To express that result in the chosen metrics, let $J_m^\pm$ be the isometric embeddings of $E_{m\pm1}$ into $\pp\otimes E_m$ supplied by these compact projections and their fixed phases. The lower map exists only for $m\geq1$.
The normalized adjacent action is
\begin{equation}\label{eq:adjacent}
R_mu_m=\left(m+\frac53\right)J_m^+u_{m+1}
+\left(m+\frac13\right)J_m^-u_{m-1}\quad(m\geq1).
\end{equation}
The coefficients are Kable's summed pairings expressed in isometric embeddings: the upper pairing is divided by $2$, and the lower by $\frac{2(m+1)}{m}$, yielding the two scalars in \eqref{eq:adjacent}. At $m=0$ only $\frac53J_0^+u_1$ remains; no lower map or singular lower coefficient is evaluated. Appendix~\ref{app:embedding-cancel} verifies the retained coefficients, the six cancellations, and the bottom action. Together, \eqref{eq:transfer} and \eqref{eq:adjacent} supply Kable's noncompact action in the normalization needed below.

\subsection{The invariant Hermitian form}
We now determine the invariant Hermitian form in the compact realization just fixed.

Use normalized Haar measure $\mathrm dk$ on the compact group $K$, with $\int_K\mathrm dk=1$, as a reference. By the Peter--Weyl (Schur) orthogonality relations, integrating the matrix coefficients in \eqref{eq:matrix-coefficient} gives
\begin{equation}\label{eq:Haar}
\langle\varphi_m(\xi),\varphi_m(\eta)\rangle_{L^2(K)}
=\frac{\|u_m\|^2}{d_m}\langle\xi,\eta\rangle_{E_m},
\quad d_m:=\dim E_m=3(m+1)^2.
\end{equation}
The $K$-types $V_m$ are pairwise inequivalent. Schur's lemma therefore makes distinct types orthogonal for any $K$-invariant Hermitian form and makes its restriction to each $V_m$ a scalar multiple of the reference form. Thus, on finite sums, it has the form
\begin{equation}\label{eq:invariant-form}
\langle v,w\rangle_{\min}
=\sum_{m\geq0}c_m\langle v_m,w_m\rangle_{L^2(K)},
\quad v_m,w_m\in V_m.
\end{equation}
The minimal representation is known to be unitary \cite[Proposition 14.2]{wallach_quaternionic_1996}, so it possesses a positive invariant form. For that form, each $c_m$ is positive. The noncompact action determines the ratios of these coefficients.

Unitarity requires, on these smooth vectors,
\begin{equation}\label{eq:adjointness}
\langle\mathrm d\pi_{\min}(Y)v,w\rangle_{\min}
=-\langle v,\mathrm d\pi_{\min}(\bar Y)w\rangle_{\min}.
\end{equation}

Set $Y=Y_{3,1}$ and
\[
v=\varphi_m(\xi^{3m+2}_{3m+2}\boxtimes\xi^m_m),\quad
w=\varphi_{m+1}(\xi^{3m+5}_{3m+5}\boxtimes\xi^{m+1}_{m+1}).
\]
The tensor of highest vectors projects to $E_{m+1}$ with coefficient one. Equations~\eqref{eq:transfer} and~\eqref{eq:adjacent} therefore give $\mathrm d\pi_{\min}(Y)v=-(m+\frac53)w$. The upward inner product is
\begin{equation}\label{eq:form-up}
\langle\mathrm d\pi_{\min}(Y)v,w\rangle_{\min}
=-\left(m+\frac53\right)c_{m+1}\frac{\|u_{m+1}\|^2}{d_{m+1}}.
\end{equation}

For the reverse direction, Kable's real structure gives $\bar Y=-Y_{-3,-1}$. The relevant lower compact projection is
\[
(J_{m+1}^-)^*\bigl(\bar Y\otimes
(\xi^{3m+5}_{3m+5}\boxtimes\xi^{m+1}_{m+1})\bigr)
=-\frac{m+1}{m+2}\,
\xi^{3m+2}_{3m+2}\boxtimes\xi^m_m.
\]
Indeed, each of the two extreme compact projections contributes $\sqrt{\frac{m+1}{m+2}}$, while conjugation supplies the minus sign. Only this component can pair with $v$; the other compact components are orthogonal. At level $m+1$ the lower action coefficient is $m+\frac43$, and the transfer minus sign cancels the sign of the projection. Consequently
\begin{equation}\label{eq:form-down}
\langle v,\mathrm d\pi_{\min}(\bar Y)w\rangle_{\min}
=\frac{m+1}{m+2}\left(m+\frac43\right)c_m\frac{\|u_m\|^2}{d_m}.
\end{equation}

Substitution of \eqref{eq:form-up} and~\eqref{eq:form-down} into adjointness gives
\[
\left(m+\frac53\right)c_{m+1}\frac{\|u_{m+1}\|^2}{3(m+2)^2}
=\frac{m+1}{m+2}\left(m+\frac43\right)c_m\frac{\|u_m\|^2}{3(m+1)^2}.
\]
Thus the ratio is
\begin{equation}\label{eq:c-ratio}
\frac{c_{m+1}}{c_m}
=\frac{(m+\frac43)(m+2)\|u_m\|^2}
{(m+\frac53)(m+1)\|u_{m+1}\|^2}.
\end{equation}
Multiplying the ratios telescopes both the dimension factors and the embedding-vector norms:
\begin{equation}\label{eq:c-product}
c_m=c_0\prod_{r=0}^{m-1}
\frac{(r+\frac43)(r+2)\|u_r\|^2}{(r+\frac53)(r+1)\|u_{r+1}\|^2}
=c_0\frac{(\frac43)_m}{(\frac53)_m}(m+1)\frac{\|u_0\|^2}{\|u_m\|^2}.
\end{equation}
Finally, fix the overall scale by agreement with the Haar form on $V_0$, that is, $c_0=1$. From \eqref{equ:4.2}, $u_0=\xi^2_0\boxtimes\xi^0_0$, and \eqref{eq:compact-norm} gives $\|u_0\|^2=\frac12$. In this normalization the coefficient in \eqref{eq:invariant-form} is therefore
\begin{equation}\label{eq:c-final}
c_m=\frac{m+1}{2\|u_m\|^2}\,
\frac{(\frac43)_m}{(\frac53)_m}.
\end{equation}
All adjacent coefficients used in the comparison are nonzero. Hence the known invariant form must have these successive ratios, and its bottom scale fixes it on every compact type. Its completion is the Hilbert space $\mathcal H_{\min}$ of the unitary minimal representation. This determines the invariant Hermitian form.

\section{The quadratic Casimir}\label{sec:basis}
\subsection{The algebraic invariant basis}
Having determined one copy of the minimal representation, we now ask for the $K$-invariant vectors in its Hilbert tensor square
\[
\mathcal H=\mathcal H_{\min}\widehat\otimes\mathcal H_{\min},
\quad \Pi(g)=\pi_{\min}(g)\otimes\pi_{\min}(g).
\]
By \eqref{eq:ladder} and compact self-duality, $(E_m\otimes E_n)^K=0$ whenever $m\ne n$. When $m=n$, this is the one-dimensional space of $K$-invariant vectors. Thus only matching compact types contribute; the orthogonal sum of these spaces is dense in $\mathcal H^K$.

In $\Gamma_q\otimes\Gamma_q$, invariance under the Cartan generator $A_1\otimes1+1\otimes A_1$ requires total weight zero. A candidate invariant vector therefore has the form
\[
I_q=\sum_{k=0}^q\alpha_k\,\xi^q_{q-2k}\otimes\xi^q_{2k-q}.
\]
Applying the raising operator $A_+\otimes1+1\otimes A_+$ and using \eqref{eq:compact-action} gives
\[
\sum_{k=0}^{q-1}\bigl((2q-2k)\alpha_k+2(k+1)\alpha_{k+1}\bigr)
\xi^q_{q-2k}\otimes\xi^q_{2k-q+2}.
\]
Here the first term in the coefficient comes from raising the second factor of the $k$-th summand, and the second from raising the first factor of the $(k+1)$-st summand. These output tensors are linearly independent, so invariance requires
\[
(2q-2k)\alpha_k+2(k+1)\alpha_{k+1}=0
\quad(0\leq k<q).
\]
Fixing $\alpha_0=1$ gives $\alpha_{k+1}=-\frac{q-k}{k+1}\alpha_k$, hence $\alpha_k=(-1)^k\binom qk$. Thus
\begin{equation}\label{eq:compact-invariant}
I_q=\sum_{k=0}^q(-1)^k\binom qk\xi^q_{q-2k}\otimes\xi^q_{2k-q}.
\end{equation}
The lowering equation gives the same coefficient relation. Thus $I_q$ is killed by the compact generators and is nonzero. Compact self-duality and multiplicity one show that it spans the one-dimensional invariant subspace. For $q=0$ the sum consists of the single trivial tensor, with no coefficient equation to impose.

To transfer this construction to the realized representation, write
\begin{equation}\label{eq:e-basis}
e^n_{k,l}=\varphi_n(\xi^{3n+2}_{3n+2-2k}\boxtimes\xi^n_{n-2l}),
\quad 0\leq k\leq3n+2,\quad 0\leq l\leq n.
\end{equation}
The product of $I_{3n+2}$ and $I_n$, with the tensor factors regrouped, gives
\begin{equation}\label{eq:Psi}
\Psi_n=\sum_{k=0}^{3n+2}\sum_{l=0}^{n}(-1)^{k+l}\binom{3n+2}{k}\binom nl
 e^n_{k,l}\otimes e^n_{3n+2-k,n-l}.
\end{equation}
The first sum runs over the weights of $\Gamma_{3n+2}$ and the second over those of $\Gamma_n$. Opposite weights are paired in each factor. The signs and binomial coefficients are exactly those required by the two raising equations; the coefficient at $k=l=0$ is one. This fixes the algebraic phase and scale of $\Psi_n$ without a Hilbert normalization.

We will calculate on the dense algebraic subspace
\begin{equation}\label{eq:algebraic-domain}
\mathcal D=(V_{\min}\otimes_{\mathrm{alg}}V_{\min})^K
=\spfin\{\Psi_n:n\geq0\}.
\end{equation}

\subsection{The quadratic contraction}\label{sec:casimir}
The next question is how the quadratic Casimir acts on \eqref{eq:Psi}. For $X\in\g$, the differential of the diagonal action is
\[
\mathrm d\Pi(X)=\mathrm d\pi_{\min}(X)\otimes1+1\otimes\mathrm d\pi_{\min}(X).
\]
To compute each ordered product in the noncompact Casimir, we apply the one-step action formulas successively. The first operator acts on the explicit sum \eqref{eq:Psi}; the second acts on every term of the resulting expression. We then collect coefficients with the same output compact types. The calculation below follows this successive-action procedure, with the repeated component expansions in Appendix~\ref{sec:app.action}.

We use the quadratic Casimir defined by $B=\frac18B_{\g}$ and a bilinear dual basis. Every compact generator annihilates a $K$-fixed vector. Hence the compact part of the Casimir vanishes on $\mathcal D$, leaving the noncompact contraction $\Cas_{\pp}$.
The dual-basis coefficients follow from the coroot Gram matrix:
\begin{equation}\label{eq:coroot-gram}
\operatorname{Gram}_B(H_\alpha,H_\beta)=\begin{pmatrix}6&-3\\-3&2\end{pmatrix},\quad
H_\vartheta=H_\alpha+2H_\beta,\quad B(H_\vartheta,H_\alpha)=0.
\end{equation}
Here $H_\vartheta=H_\alpha+2H_\beta$ is the highest-root coroot. In the weight basis $Y_{a,b}$ defined before \eqref{eq:R-zeta}, the isometry \eqref{eq:iota} gives
\[
\|Y_{3,1}\|^2=\|Y_{3,-1}\|^2=1,
\quad \|Y_{1,1}\|^2=\|Y_{1,-1}\|^2=\frac13.
\]
The conjugate weights have the same norms, and distinct compact weights are Hermitian orthogonal. In Kable's real structure, the conjugates of $Y_{3,1},Y_{3,-1},Y_{1,1},Y_{1,-1}$ are respectively $-Y_{-3,-1},Y_{-3,1},Y_{-1,-1},-Y_{-1,1}$. These signs and norms determine the bilinear duals.
The $B$-dual of $Y_{a,b}$ is $\frac{\overline{Y_{a,b}}}{\|Y_{a,b}\|^2}$, giving
\begin{align}
\Cas_{\pp}&=-Y_{3,1}Y_{-3,-1}+Y_{3,-1}Y_{-3,1}
+3Y_{1,1}Y_{-1,-1}-3Y_{1,-1}Y_{-1,1}\nonumber\\
&\quad-Y_{-3,-1}Y_{3,1}+Y_{-3,1}Y_{3,-1}
+3Y_{-1,-1}Y_{1,1}-3Y_{-1,1}Y_{1,-1}.\label{equ:4.cas_expansion}
\end{align}

The adjacent action \eqref{eq:adjacent} shifts a single factor by one level. Two actions on one factor can return to its original level or shift it by two; mixed actions shift each factor by one. Since the Casimir preserves compact invariance, only matching levels survive. These are the original pair and its two neighbors. This explains the tridiagonal form; the following theorem gives its coefficients.

\begin{theorem}[Diagonal Casimir recurrence]\label{thm:4.cas_action}
In the algebraic invariant basis \eqref{eq:Psi}, the quadratic Casimir for the diagonal action satisfies, for $n\geq0$,
\begin{equation}\label{equ:4.recurrence_final}
\mathrm d\Pi(\Cas)\Psi_n=-2(n+\frac{5}{3})^2\Psi_{n+1}-2\delta_n\Psi_n-2(n+\frac{1}{3})^2\Psi_{n-1},
\end{equation}
where $\Psi_{-1}=0$ and
\begin{equation}\label{eq:delta}
\delta_n=\frac{(n+2)(n+\frac{4}{3})(n+\frac{5}{3})+n(n+\frac{1}{3})(n+\frac{2}{3})}{n+1}.
\end{equation}
\end{theorem}
\begin{proof}
For $Y,Z\in\pp$, the difference of the two ordered actions on $\Psi_n$ is
\[
\bigl(\mathrm d\Pi(Y)\mathrm d\Pi(Z)-\mathrm d\Pi(Z)\mathrm d\Pi(Y)\bigr)\Psi_n
=\mathrm d\Pi([Y,Z])\Psi_n=0,
\]
since $[\pp,\pp]\subset\kk$ and $\Psi_n$ is $K$-invariant. Each reversed product in \eqref{equ:4.cas_expansion} therefore equals its forward product on this vector. It suffices to calculate the four signed representative products
\begin{equation}\label{eq:Q}
Q=-\mathcal Y_{3,1}\mathcal Y_{-3,-1}+\mathcal Y_{3,-1}\mathcal Y_{-3,1}+3\mathcal Y_{1,1}\mathcal Y_{-1,-1}-3\mathcal Y_{1,-1}\mathcal Y_{-1,1},
\end{equation}
where $\mathcal Y_{a,b}=\mathrm d\Pi(Y_{a,b})$ and as a consequence $\frac12\mathrm d\Pi(\Cas)\Psi_n=Q\Psi_n$.

For each summand
\[
e^n_{k,l}\otimes e^n_{3n+2-k,n-l}
\]
of $\Psi_n$, apply the first operator of each ordered product in $Q$ using the adjacent one-step action, and then apply the second operator to every term produced. Multiply by the binomial coefficient and sign in \eqref{eq:Psi}, sum over $k,l$, and group the outputs by their pair of compact levels. Appendix~\ref{app:paths} gives the first-order formulas, an explicit successive-action example, and all second-order paths.

The matching output pairs are $(n+1,n+1)$, $(n,n)$, and $(n-1,n-1)$. The remaining pairs are $(n\pm2,n)$, $(n,n\pm2)$, $(n+1,n-1)$, and $(n-1,n+1)$, whenever those levels exist. Appendix~\ref{app:nonmatching} verifies that every nonmatching group cancels after the four signed products are combined. Thus these groups are canceled in the computation, rather than omitted from the intermediate actions.

The three matching groups give the coefficients of $\Psi_{n+1}$, $\Psi_n$, and $\Psi_{n-1}$. Appendix~\ref{app:aggregation} performs their reindexing and aggregation. Both factors move upward in the $\Psi_{n+1}$ group, giving $-(n+\frac53)^2$, and both move downward in the $\Psi_{n-1}$ group, giving $-(n+\frac13)^2$. A return to $(n,n)$ comes from two successive actions on the same factor, either up--down or down--up. Including both choices of tensor factor, the negatives of the two diagonal contributions are respectively
\[
\frac{(n+2)(n+\frac43)(n+\frac53)}{n+1}
\quad\text{and}\quad
\frac{n(n+\frac13)(n+\frac23)}{n+1}.
\]
Their sum is $\delta_n$. The down--up group is absent at $n=0$.

Combining the three retained groups gives
\[
Q\Psi_n=-(n+\tfrac53)^2\Psi_{n+1}-\delta_n\Psi_n-(n+\tfrac13)^2\Psi_{n-1}.
\]
By the reversed-order equality at the start of the proof, the complete Casimir has twice these coefficients. This proves \eqref{equ:4.recurrence_final} at generic levels.

At $n=0$ and $n=1$ some compact summands do not exist. They must be omitted before evaluating any coefficient with a vanishing denominator. At zero, the first action has only upward terms, and the second action from level one is regular. At one, a first downward step reaches zero and cannot be followed by another downward step in that factor. Applying only the existing one-step formulas gives the four signed products listed in Appendix~\ref{app:boundary}; their sums, doubled as in \eqref{eq:Q}, are
\begin{align}\label{eq:casimir-boundary}
&\mathrm d\Pi(\Cas)\Psi_0=-\tfrac{50}{9}\Psi_1-\tfrac{80}{9}\Psi_0,\nonumber\\
&\mathrm d\Pi(\Cas)\Psi_1=-\tfrac{128}{9}\Psi_2-\tfrac{188}{9}\Psi_1-\tfrac{32}{9}\Psi_0.
\end{align}
These regular calculations agree with the same recurrence and complete the proof.
\end{proof}

\section{Wilson spectral realization}\label{sec:spectral}
The algebraic recurrence must now be expressed in an orthonormal basis. For the vector $e^n_{k,l}$ in \eqref{eq:e-basis}, the invariant form \eqref{eq:invariant-form} gives
\[
\|e^n_{k,l}\|_{\min}^2=c_n\|e^n_{k,l}\|_{L^2(K)}^2.
\]
By the Peter--Weyl formula \eqref{eq:Haar},
\[
\|e^n_{k,l}\|_{L^2(K)}^2
=\frac{\|u_n\|^2}{\dim E_n}
\|\xi^{3n+2}_{3n+2-2k}\boxtimes\xi^n_{n-2l}\|^2.
\]
We therefore put
\[
h_n:=\frac{c_n\|u_n\|^2}{\dim E_n}
=\frac{c_n\|u_n\|^2}{3(n+1)^2}.
\]
This abbreviation is introduced only to avoid repeating the common scalar factor in the norms of all vectors in the $n$-th $K$-type. Substituting \eqref{eq:c-final} and using the compact weight norms \eqref{eq:compact-norm} gives
\begin{equation}\label{eq:h}
h_n=\frac{(\frac43)_n}{6(n+1)(\frac53)_n},\quad
\|e^n_{k,l}\|_{\min}^2=\frac{h_n}{\binom{3n+2}{k}\binom nl}.
\end{equation}
The tensors in \eqref{eq:Psi} are mutually orthogonal. Opposite weights have equal norm, since $\binom{3n+2}{3n+2-k}=\binom{3n+2}{k}$ and $\binom n{n-l}=\binom nl$. Each summand therefore contributes
\[
\binom{3n+2}{k}^{\!2}\binom nl^{\!2}\!
\|e^n_{k,l}\|_{\min}^2\|e^n_{3n+2-k,n-l}\|_{\min}^2\!=\!\binom{3n+2}{k}^{\!2}\binom nl^{\!2}\!
\frac{h_n}{\binom{3n+2}{k}\binom nl}
\frac{h_n}{\binom{3n+2}{k}\binom nl}\!=\!h_n^2.
\]
There are $(3n+3)(n+1)=3(n+1)^2$ terms, whence
\[
\|\Psi_n\|^2=\sum_{k=0}^{3n+2}\sum_{l=0}^n h_n^2=3(n+1)^2h_n^2.
\]
Taking the positive square root gives
\begin{equation}\label{eq:rho}
\rho_n=\|\Psi_n\|=\sqrt3(n+1)h_n
=\frac{(\frac43)_n}{2\sqrt3(\frac53)_n}.
\end{equation}
We now define $\Phi_n=\frac{\Psi_n}{\rho_n}$.
Thus
\begin{equation}\label{eq:HK-domain}
\mathcal H^K=\widehat\bigoplus_{n\geq0}\C\Phi_n,
\quad \mathcal D=\spfin\{\Phi_n:n\geq0\}.
\end{equation}
This normalization turns the tridiagonal Casimir matrix into a symmetric Jacobi matrix. We state the resulting spectral realization first, then prove it by this symmetrization and a comparison with the Wilson recurrence.

\begin{theorem}[Wilson realization of the diagonal Casimir]\label{thm:4.spectral_decomp}
Write $\mathcal H^\infty$ for the smooth vectors of the diagonal representation on $\mathcal H$. Let
\[
\mathcal C=\overline{\mathrm d\Pi(\Cas)|_{\mathcal D}}
=\overline{\mathrm d\Pi(\Cas)|_{\mathcal H^\infty\cap \mathcal H^K}}.
\]
This is a self-adjoint operator on $\mathcal H^K$. Set
\begin{equation}\label{eq:Wilson-parameters}
(a,b,c,d)=\left(\frac34,\frac34,\frac7{12},\frac{11}{12}\right).
\end{equation}
For $x>0$, define the normalized spectral density
\begin{equation}\label{eq:density}
\omega(x)=\frac1\pi\left|\frac{\Gamma(\frac34+ix)^2\Gamma(\frac7{12}+ix)\Gamma(\frac{11}{12}+ix)}{\Gamma(2ix)\Gamma(\frac32)\Gamma(\frac43)\Gamma(\frac53)}\right|^2.
\end{equation}
Then the assignment
\begin{equation}\label{eq:transform}
T\Phi_n=(n+1)(-1)^n\widetilde W_n(x^2;a,b,c,d)
\end{equation}
extends to a unitary isomorphism $T:\mathcal H^K\to L^2((0,\infty),\omega(x)\dd x)$, and
\begin{equation}\label{eq:operator-identity}
T\mathcal C T^{-1}=M_{-\frac{9}{2}-8x^2}.
\end{equation}
Here $M_f$ denotes multiplication by the real function $f$. The multiplication operator has its maximal domain. Explicitly,
\begin{equation}\label{eq:Cas-domain}
\Dom\mathcal C=\biggl\{f\in \mathcal H^K:\biggl(\frac{9}{2}+8x^2\biggr)Tf\in L^2((0,\infty),\omega(x)\dd x)\biggr\}.
\end{equation}
The finite span $\mathcal D$ is a core.
\end{theorem}

The proof proceeds in the following two steps.
\subsection{Symmetrization and self-adjoint realization}
On the finite domain $\mathcal D$ in \eqref{eq:HK-domain}, define
\[
A_0=-\tfrac12\mathrm d\Pi(\Cas)|_{\mathcal D}.
\]
Using \eqref{eq:rho}, the recurrence becomes
\begin{equation}\label{equ:4.Phi_recurrence}
A_0\Phi_n=a_n\Phi_{n+1}+\delta_n\Phi_n+a_{n-1}\Phi_{n-1},\quad
a_n=(n+\frac{5}{3})(n+\frac{4}{3}),
\end{equation}
with the last term omitted at $n=0$. This follows from
\[
\begin{aligned}
&(n+\frac{5}{3})^2\frac{\rho_{n+1}}{\rho_n}=(n+\frac{5}{3})(n+\frac{4}{3}),\\
&(n+\frac{1}{3})^2\frac{\rho_{n-1}}{\rho_n}=(n+\frac{1}{3})(n+\frac{2}{3})\quad(n\geq1).
\end{aligned}
\]

\begin{lemma}\label{lem:closure}
The finite-domain operator $A_0$ is essentially self-adjoint. Its closure equals the closure of $S=-\tfrac12\mathrm d\Pi(\Cas)$ on $\mathcal H^\infty\cap \mathcal H^K$, where $\mathcal H^\infty$ denotes smooth vectors for the diagonal representation. In particular $\mathcal D$ is a core for this closure.
\end{lemma}
\begin{proof}
Under $\Phi_n\leftrightarrow\delta_n^{\mathrm{basis}}$ in $\ell^2(\N_0)$, the Jacobi expression is
\[
(\tau f)_n=a_nf_{n+1}+\delta_nf_n+a_{n-1}f_{n-1},
\]
with a fixed left endpoint and no last term at zero. The minimal Jacobi operator is the closure of this expression on finite-support sequences. The maximal operator has domain $\{f\in\ell^2:\tau f\in\ell^2\}$. The coefficients are real and $a_n>0$. The explicit sequence $v_n=(-1)^n(n+1)$ solves $\tau v=0$, since
\begin{equation}\label{eq:zero-solution}
\delta_n(n+1)=a_n(n+2)+a_{n-1}n\quad(n\geq1),\quad\delta_0=\frac{40}{9}=2a_0.
\end{equation}
It is not square summable. The half-line limit-point criterion therefore makes the maximal operator self-adjoint and equal to the closure of the finite-support operator; see \cite[Section 2.6, Lemmas 2.15--2.16, pp.\,48--49]{teschl_jacobi_2000}, including the minimal/maximal adjoint relations on pp.\,47--48. The criterion applies at a real spectral parameter.

Now let $S=-\frac12\mathrm d\Pi(\Cas)$ on $\mathcal H^\infty\cap\mathcal H^K$. The real Casimir is symmetric in the differential unitary representation. Centrality preserves $\mathcal H^K$, and finite tensors of Harish-Chandra vectors are smooth. Thus $A_0\subset S$ and symmetry gives
\[
A_0\subset S\subset A_0^*=\overline{A_0}.
\]
Taking closures proves $\overline S=\overline{A_0}$ and establishes the finite algebraic span as a core.
\end{proof}

\subsection{Recurrence matching and the spectral transform}
We compare the symmetric recurrence directly with the polynomial images of the orthonormal vectors $\Phi_n$. Define
\begin{equation}\label{eq:p-polynomials}
p_n(x^2)=(n+1)(-1)^n\widetilde W_n(x^2;a,b,c,d),
\end{equation}
with the parameters in \eqref{eq:Wilson-parameters}. The parameter combinations in the general Wilson coefficients are $
a+b+c+d=\frac34+\frac34+\frac7{12}+\frac{11}{12}=3$,
\[
a+b=c+d=\frac32,\quad a+c=b+c=\frac43,\quad a+d=b+d=\frac53.
\]
Substitution in \eqref{eq:Wilson-coefficients}, followed by cancellation of the factors $n+\frac32$ and $n+\frac12$, gives, for $n\geq1$,
\begin{equation}\label{eq:Wilson-specialized}
\begin{aligned}
&A_n^W=\frac{(n+2)(n+\frac32)(n+\frac43)(n+\frac53)}{(2n+2)(2n+3)}
=\frac{(n+2)(n+\frac43)(n+\frac53)}{4(n+1)},\\
&C_n^W=\frac{n(n+\frac13)(n+\frac23)(n+\frac12)}{(2n+1)(2n+2)}
=\frac{n(n+\frac13)(n+\frac23)}{4(n+1)}.
\end{aligned}
\end{equation}
Multiplying \eqref{eq:Wilson-recurrence} by $-4(n+1)(-1)^n$ and expressing each term in \eqref{eq:p-polynomials} yields
\begin{equation}\label{eq:p-recurrence}
4(a^2+x^2)p_n=4A_n^W\frac{n+1}{n+2}p_{n+1}
+4(A_n^W+C_n^W)p_n+4C_n^W\frac{n+1}{n}p_{n-1},\quad n\geq1.
\end{equation}
Each coefficient now agrees with the symmetric Jacobi matrix:
\[
\begin{aligned}
&4A_n^W\frac{n+1}{n+2}=(n+\tfrac53)(n+\tfrac43)=a_n,\\
&4C_n^W\frac{n+1}{n}=(n+\tfrac13)(n+\tfrac23)=a_{n-1},\\
&4(A_n^W+C_n^W)=
\frac{(n+2)(n+\frac43)(n+\frac53)+n(n+\frac13)(n+\frac23)}{n+1}
=\delta_n.
\end{aligned}
\]
At zero, \eqref{eq:Wilson-initial} specializes to $A_0^W=\frac{10}{9}$ and $p_1=-2\widetilde W_1$, so
\begin{equation}\label{eq:Wilson-zero}
4(a^2+x^2)p_0=\frac{20}{9}p_1+\frac{40}{9}p_0=a_0p_1+\delta_0p_0.
\end{equation}
There is no lower term and no evaluation of the quotient containing $n$ in \eqref{eq:p-recurrence} at zero.

The two normalizations have different sources. The passage $\Psi_n\mapsto\Phi_n$ uses the representation's Hilbert norm and symmetrizes the Casimir matrix. The factor $n+1$ in $p_n$ adjusts the Wilson polynomial normalization to that same matrix; its norm will be checked below. Since $4(a^2+x^2)=\frac94+4x^2$, the Casimir itself has multiplier $-\frac92-8x^2$.

The recurrence has identified a multiplication model; its measure must now be normalized so that the images of $\Phi_n$ form an orthonormal basis. Write
\[
C_\Gamma=\Gamma(\tfrac32)\Gamma(\tfrac43)\Gamma(\tfrac53).
\]
At these parameters, the six pairwise sums in Wilson's squared-norm formula are two copies each of $\frac32,\frac43,\frac53$. Using the weight \eqref{eq:Wilson-weight} and canceling the polynomial normalization factors gives
\begin{equation}\label{eq:Wilson-norm}
\|\widetilde W_n\|^2_{L^2(\omega_W)}
=C_\Gamma^2\frac{\Gamma(2n+2)\Gamma(n+1)}{\Gamma(n+2)\Gamma(2n+3)}
=\frac{C_\Gamma^2}{2(n+1)^2}.
\end{equation}
It follows that $\omega=\frac{2\omega_W}{C_\Gamma^2}$ is precisely \eqref{eq:density}. Consequently, by \eqref{eq:p-polynomials} and \eqref{eq:Wilson-norm},
\[
\|p_n\|_{L^2(\omega)}^2
=(n+1)^2\frac{2}{C_\Gamma^2}
\frac{C_\Gamma^2}{2(n+1)^2}=1.
\]
Lemma~\ref{lem:Wilson-complete} supplies both orthogonality and completeness in the full $L^2$ space. Neither property changes when the measure is multiplied by a positive constant or the individual polynomials by nonzero constants.

\begin{proof}[Completion of the proof of Theorem $\ref{thm:4.spectral_decomp}$]
The assignment \eqref{eq:transform} maps the orthonormal basis of $\mathcal H^K$ onto the complete orthonormal polynomial system just obtained. It thus extends uniquely to a unitary map onto the full weighted $L^2$ space. In particular $T\Phi_0=1$ and $\int_0^\infty\omega(x)\dd x=1$.

The recurrence comparison, including \eqref{eq:Wilson-zero}, intertwines $A_0$ on $\mathcal D$ with multiplication by $\frac94+4x^2$ on the polynomial span. Maximal multiplication by this real function is self-adjoint and extends that polynomial-domain operator. Lemma~\ref{lem:closure} gives uniqueness of the self-adjoint extension, hence
\[
T\overline{A_0}T^{-1}=M_{\frac94+4x^2}.
\]
Multiplying by $-2$ proves \eqref{eq:operator-identity}. The maximal multiplication domain is exactly the domain stated in \eqref{eq:Cas-domain}. The equality with the smooth representation-theoretic closure and the core assertion follow from Lemma~\ref{lem:closure}.
\end{proof}

\begin{corollary}[Spectral consequence]\label{cor:spectrum}
The spectrum of $\mathcal C$ on $\mathcal H^K$ is absolutely continuous, equals $(-\infty,-\frac{9}{2}]$, and has multiplicity one. There are no discrete atoms.
\end{corollary}
\begin{proof}
The density is positive for $x>0$, and $\lambda(x)=-\frac{9}{2}-8x^2$ is one-to-one there. Its closed essential range is the stated interval. The push-forward density on the open interval is $\frac{\omega(x)}{16x}$. The endpoint has no singleton mass. The scalar multiplication model has spectral multiplicity one; $T\Phi_0=1$ and polynomial completeness also give a cyclic vector.
\end{proof}

\section{A conjectural representation interpretation}\label{sec:conjectural}
The Wilson transform of Theorem~\ref{thm:4.spectral_decomp} describes the quadratic Casimir on
\[
\mathcal H^K=(\mathcal H_{\min}\widehat\otimes\mathcal H_{\min})^K.
\]
It suggests a representation-theoretic interpretation of the closed $G$-subspace generated by these vectors,
\[
\mathcal H_{\mathrm{sph}}
:=\overline{\operatorname{span}\bigl(\Pi(G)\mathcal H^K\bigr)}
\subset\mathcal H_{\min}\widehat\otimes\mathcal H_{\min}.
\]

Write the Heisenberg parabolic $Q=LN$ in Langlands form as $P=MAN$, so that $P=Q$ and $L=MA$. Following the notation for the degenerate principal series in \cite{frahm_degenerate_2026}, consider the induced family
\[
\pi_{s,\varepsilon}
=\operatorname{Ind}_{P}^{G}
\bigl(\operatorname{sgn}^{\varepsilon}\otimes\chi_s\otimes1\bigr),
\quad s\in\C,\quad\varepsilon\in\{0,1\}.
\]
Here $\operatorname{sgn}$ is the nontrivial component character of $M$ and $\chi_s$ is the character of $A$ with continuous parameter $s$. We use the normalization in which the unitary principal axis is $s\in\frac32+i\R$. Only the branch $\varepsilon=0$ contains a $K$-fixed vector.

For $\pi_{s,0}$ the infinitesimal character in these conventions is
$\Lambda(s)=\rho-s\varpi_\beta$, where $\rho$ is the half-sum of the positive roots and $\varpi_\beta=3\alpha+2\beta$. With the induced bilinear form on weights normalized by
$B(\varpi_\beta,\varpi_\beta)=2$ and $B(\rho,\varpi_\beta)=3$, the quadratic Casimir scalar is
\[
B(\rho-s\varpi_\beta,\rho-s\varpi_\beta)-B(\rho,\rho)
=2s^2-6s.
\]
Writing $s=\frac32+i\lambda$ gives
\[
2s^2-6s=-\frac92-2\lambda^2.
\]
Comparison with the proved Wilson multiplier $-\frac92-8x^2$ therefore gives $\lambda=2x$ and
\[
s=\frac32+2ix,
\]
with the conjugate choice $s=\frac32-2ix$ giving the same scalar. The positive parameter $x$ is compatible with this symmetry and with the support of the Wilson measure. This comparison motivates the following conjecture.

\begin{conjecture}\label{conj:4.spherical_plancherel}
The closed spherical-generated representation $\mathcal H_{\mathrm{sph}}$ is multiplicity free and admits a $G$-equivariant unitary direct-integral decomposition
\[
\mathcal H_{\mathrm{sph}}
\simeq\int_0^\infty{}^{\oplus}
\pi_{\frac32+2ix,\,0}\,\omega(x)\dd x,
\]
where $\omega(x)\dd x$ is the normalized Wilson measure in Theorem $\ref{thm:4.spectral_decomp}$.
\end{conjecture}

The conjectured decomposition would induce the proved scalar model on $K$-fixed vectors; the converse has not been established. Besides the quadratic-Casimir match, the Gamma-factor form of the Wilson density resembles a spherical $c$-function density, and its absence of atoms is compatible with the proposed continuous direct integral. These are supporting features, not a construction of a $G$-intertwining spherical Fourier transform.

For $G_2$, the quadratic Casimir alone does not determine the full infinitesimal character. What remains is to construct or identify a $G$-equivariant direct-integral transform realizing these fibers, or otherwise determine $\mathcal H_{\mathrm{sph}}$. Further central data, such as the independent degree-six central element, or a direct intertwining construction provide natural approaches. We leave this representation-theoretic identification as an open problem.

\appendix
\section{Compact normalization and adjacent action}\label{app:compact}
We record the projections of the maps $f_m,g_n$ in \eqref{eq:CG-interface}, then their normalization and their application to $R_mu_m$. For an isometric decomposition, changing from orthonormal weight vectors to \eqref{eq:compact-norm} multiplies a coefficient into $\Gamma_h$ by
\[
\left(\frac{\binom hj}{\binom nk\binom qr}\right)^{\frac{1}{2}},
\quad j=k+r-\frac{n+q-h}{2},
\]
for the input $\xi^n_{n-2k}\otimes\xi^q_{q-2r}$. Only existing summands and valid weights are evaluated.

The general orthonormal Clebsch--Gordan coefficients used in this conversion are those of \cite[Section 8.1.3]{vilenkin_representation_2012}; their spins are $\frac n2,\frac q2,\frac h2$. The following two normalizations specify the phases and scales needed here.

\subsection{The \texorpdfstring{$\Gamma_1$}{Gamma 1} normalization}
Flip the two factors of $\Gamma_p\otimes\Gamma_q$ as necessary to place $\pp$ first in $\pp\otimes E_n$; the flip is the ordinary unsigned tensor flip. For $\Gamma_m\otimes\Gamma_1$, the projection onto the upper summand has coefficient one. The lower coefficients are
\[
\begin{aligned}
&\xi^m_{m-2l}\otimes\xi^1_1\longmapsto-\frac{l}{\sqrt{m(m+1)}}\xi^{m-1}_{m+1-2l},\\
&\xi^m_{m-2l}\otimes\xi^1_{-1}\longmapsto\frac{m-l}{\sqrt{m(m+1)}}\xi^{m-1}_{m-1-2l}.\end{aligned}
\]
for $m\geq1$. For example, the orthonormal coefficient for $\xi^m_{m-2l}\otimes\xi^1_{-1}$ into $\Gamma_{m-1}$ is $\sqrt{\frac{m-l}{m+1}}$. The change of basis contributes $\sqrt{\frac{m-l}{m}}$, giving $\frac{m-l}{\sqrt{m(m+1)}}$. The sign for the other weight follows from the raising relation. This checks the $\Gamma_1$ normalization directly.

\subsection{The \texorpdfstring{$\Gamma_3$}{Gamma 3} normalization}
For $\Gamma_n\otimes\Gamma_3$, the projection onto $\Gamma_{n+3}$ again has coefficient one. We use the falling factorial $(k)^{\underline\ell}=k(k-1)\cdots(k-\ell+1)$ for nonnegative integer $\ell$, with $(k)^{\underline0}=1$. Writing $a=3-2r$, the extreme lower projection is
\begin{equation}\label{eq:CG-extreme}
\xi^n_{n-2k}\otimes\xi^3_{3-2r}\longmapsto
\frac{(-1)^{r+1}(k)^{\underline{3-r}}(n-k)^{\underline r}}
{n(n-1)\sqrt{(n+1)(n-2)}}\xi^{n-3}_{n+3-2r-2k},\quad n\geq3.
\end{equation}
The normalization of an intermediate projection can also be checked without an entire coefficient table. For $n\geq2$, the orthonormal coefficient of $\xi^n_{n-2k}\otimes\xi^3_{-1}$ into $\Gamma_{n-1}$ is
\[
(n-3k-1)\sqrt{\frac{n-k}{(n-1)(n+1)(n+2)}}.
\]
Changing bases contributes $\sqrt{\frac{n-k}{3n}}$, and the resulting coefficient is
\[
\frac{(n-k)(n-3k-1)}{\sqrt{3(n+2)(n+1)n(n-1)}}.
\]
This is the second compact normalization used in the eight projections below. The formula is used on its valid weight range; absent weights are zero.

These are the phases of the isometric projections $(J_m^\pm)^*$. Multiplying the two lower coefficients by $m+\frac{1}{3}$ and using $n=3m+2$ gives exactly \eqref{eq:L}; the transfer formula supplies the overall minus in \eqref{eq:one-factor}. At $m=0$ both lower representation factors are absent. We use only the upper action there.

\subsection{The eight projections}\label{app:embedding-cancel}
Apply the compact maps just normalized to the six terms of \eqref{eq:R-zeta}. Write $M=3m+2$ and let $w_{m,k}$ denote the coefficient of $\zeta_{m-2k}$ in \eqref{equ:4.2}, with $w_{m,k}=0$ outside $0\leq k\leq m$.
For $m\geq1$, the candidate projections are indexed by $(r,h)$, with target $\Gamma_{M+3-2r}\boxtimes\Gamma_{m+1-2h}$. Put
\begin{align*}
\mathsf D_0&=1,&\mathsf D_1&=\sqrt{3M(M+3)},\\
\mathsf D_2&=\sqrt{3(M+2)(M+1)M(M-1)},&
\mathsf D_3&=M(M-1)\sqrt{(M+1)(M-2)},\\
\mathsf E_0&=1,&\mathsf E_1&=\sqrt{m(m+1)}.
\end{align*}
Each projection of \eqref{eq:R-zeta} has weights $(3m+3-6k,m+1-2k)$ and $(3m-3-6k,m-1-2k)$, with coefficients $A_{r,h}(k)/(\mathsf D_r\mathsf E_h)$ and $B_{r,h}(k)/(\mathsf D_r\mathsf E_h)$. For $v=m-k$ the numerators are
{\small
\[
\begin{array}{c|l|l}
(r,h)&A_{r,h}(k)&B_{r,h}(k)\\\hline
(0,0)&m+2k+\frac{5}{3}&-(3m-2k+\frac{5}{3})\\
(0,1)&-\tfrac23 k(3k+1)&-\tfrac23 v(3v+1)\\
(1,0)&-6k(3k+1)&-6v(3v+1)\\
(1,1)&-3k(3k+1)(m+1-2k)&-3v(3v+1)(m-1-2k)\\
(2,0)&-9k(3k+1)(m+1-2k)&-9v(3v+1)(m-1-2k)\\
(2,1)&18k^2(3k+1)(v+1)&18v^2(3v+1)(k+1)\\
(3,0)&k(3k+1)(3k-1)(6v+4)&v(3v+1)(3v-1)(6k+4)\\
(3,1)&-k^2(3k+1)(3k-1)(9m+7-6k)&v^2(3v+1)(3v-1)(3m+6k+7).
\end{array}
\]
}
At a common output weight the coefficient in $Ru_m$ is
\[
\frac{w_{m,k}A_{r,h}(k)+w_{m,k-1}B_{r,h}(k-1)}{\mathsf D_r\mathsf E_h}.
\]
\subsection{The retained upper and lower coefficients}
First collect the upper row $(r,h)=(0,0)$ against $w_{m,k}$ and $w_{m,k-1}$. Dividing by the coefficient of the same weight in $u_{m+1}$ reduces the sum to the first identity below. Next collect the lower row $(3,1)$ against $w_{m,k}$ and $w_{m,k+1}$ and divide by the corresponding coefficient in $u_{m-1}$; this gives the second identity, for $m\geq1$:
\begin{align*}
&(m-k+1)(m-k+\frac{4}{3})(m+2k+\frac{5}{3})+k(k+\frac{1}{3})(3m-2k+\frac{11}{3})\\
&\quad =(m+1)(m+\frac{4}{3})(m+\frac{5}{3}),\\
&\frac{(m-k)(3m-3k-1)(3m+6k+7)}{(3m+3)(3m+2)}+\frac{(k+1)(3k+2)(9m-6k+1)}{(3m+3)(3m+2)}\\
&\quad=m+\frac{1}{3}.
\end{align*}
They give the two coefficients in \eqref{eq:adjacent}.

In the proof of \cite[Theorem 4.8, p.\,122]{kable_conformally_2012}, the upper test vector pairs with $J_m^+u_{m+1}$ by $2$, and its lower test vector pairs with $J_m^-u_{m-1}$ by $\frac{2(m+1)}m$. The corresponding pairings of $Ru_m$ are $2(m+\frac{5}{3})$ and $\frac{2(m+\frac{1}{3})(m+1)}m$, respectively. Division gives \eqref{eq:adjacent}. The lower CG highest vector there has norm $\frac{m+1}m$; this conversion retains its positive extreme phase.

\subsection{The six cancellations and the bottom level}
The ratio of the embedding coefficients is
\begin{equation}\label{eq:w-ratio}
\frac{w_{m,k-1}}{w_{m,k}}=-\frac{k(3k+1)}{(m-k+1)(3m-3k+4)},\quad1\leq k\leq m.
\end{equation}
For each of the six rows $(0,1),(1,0),(1,1),(2,0),(2,1),(3,0)$ the polynomial identity is
\begin{equation}\label{eq:B2-cancel}
(m-k+1)(3m-3k+4)A_{r,h}(k)-k(3k+1)B_{r,h}(k-1)=0.
\end{equation}
Direct substitution of the table proves every cancellation. Endpoint terms vanish by $A(0)=0$, $B(m)=0$, or absent weights. \paragraph{Bottom action.}
At $m=0$, omit the $h=1$ and $r=3$ summands before evaluating coefficients. The actual compact product is $(\Gamma_5\oplus\Gamma_3\oplus\Gamma_1)\boxtimes\Gamma_1$. The six-term action reduces to
\[
Ru_0=\tfrac23Y_{3,1}\otimes u_0-\tfrac23Y_{-3,-1}\otimes u_0
+Y_{1,1}\otimes(\xi_2^2\boxtimes1)-Y_{-1,-1}\otimes(\xi_{-2}^2\boxtimes1).
\]
The upper projection is $\tfrac53(\xi_3^5\boxtimes\xi_1^1-\xi_{-3}^5\boxtimes\xi_{-1}^1)=\tfrac53u_1$. The $\Gamma_3\boxtimes\Gamma_1$ terms cancel with coefficients $\pm2/\sqrt{30}$, and $\Gamma_1$ cannot carry first-factor weights $\pm3$. This verifies the regular bottom action.

\raggedbottom
\section{Tensor-Casimir coefficient calculation}\label{sec:app.action}
This appendix verifies the tensor contraction in Theorem~\ref{thm:4.cas_action}. The first two subsections give the action paths and the explicit nonmatching-sector cancellations; the third aggregates the three retained coefficients; the last checks the boundary indices. On a first reading, the long index sums may be skipped after noting the mixed factor \eqref{eq:mixed-factor} and its cancellation \eqref{eq:mixed-cancel}.

Throughout, we use
\[
U_n=n+\frac53\quad\text{and}\quad D_n=9(3n+2)n(n+1)\quad(n\geq1).
\]
These abbreviations are used only to shorten the repeated first-order action formulas. The denominator $D_n$ is used only for existing downward steps.

The first-order expressions and the four second-order blocks are generated by the following single rule, retaining all paths required for their contraction.

\subsection{First-order action and second-order paths}\label{app:paths}
For $a=3-2r$, $b=1-2t$, with $0\leq r\leq3$, $0\leq t\leq1$, define, only when $n\geq1$,
\begin{equation}\label{eq:L}
L_{r,t}(n,k,l)=\frac{(-1)^{r+t}(k)^{\underline{3-r}}(3n+2-k)^{\underline r}
(l)^{\underline{1-t}}(n-l)^{\underline t}}{D_n}.
\end{equation}
The credited adjacent action, the transfer sign, and the extreme CG maps in Appendix~\ref{app:compact} give
\begin{equation}\label{eq:one-factor}
\mathrm d\pi_{\min}(Y_{3-2r,1-2t})e^n_{k,l}
=-U_ne^{n+1}_{k+r,l+t}-L_{r,t}(n,k,l)e^{n-1}_{k+r-3,l+t-1}.
\end{equation}
At zero the lower term is omitted entirely; $L_{r,t}(0,k,l)$ is never evaluated. Invalid target weights in existing representations give zero. For clarity the eight numerators of $L$ are
\[
\begin{array}{c|l@{\qquad}c|l}
(r,t)&D_nL_{r,t}&(r,t)&D_nL_{r,t}\\\hline
(0,0)&(k)^{\underline3}l&(0,1)&-(k)^{\underline3}(n-l)\\
(1,0)&-(k)^{\underline2}(3n+2-k)l&(1,1)&(k)^{\underline2}(3n+2-k)(n-l)\\
(2,0)&k(3n+2-k)^{\underline2}l&(2,1)&-k(3n+2-k)^{\underline2}(n-l)\\
(3,0)&-(3n+2-k)^{\underline3}l&(3,1)&(3n+2-k)^{\underline3}(n-l).
\end{array}
\]

Start with a tensor $e^n_{k,l}\otimes e^n_{3n+2-k,n-l}$ occurring in \eqref{eq:Psi}. The one-step formula \eqref{eq:one-factor} gives only levels $n+1$ and $n-1$ in the factor on which the operator acts. To see the full calculation for the first ordered product, use the two coefficients
\[
\begin{aligned}
&L_{0,0}(n,k,l)=\frac{(k)^{\underline3}l}{9(3n+2)n(n+1)},\\
&L_{3,1}(n,k,l)=\frac{(3n+2-k)^{\underline3}(n-l)}{9(3n+2)n(n+1)}
\quad(n\geq1).
\end{aligned}
\]
Here $(k)^{\underline\ell}=k(k-1)\cdots(k-\ell+1)$, with $(k)^{\underline0}=1$, is the falling factorial. These are two instances of \eqref{eq:L}. The complete first action on a typical tensor is
\begin{equation}\label{eq:first-tensor-example}
\begin{aligned}
&\mathcal Y_{-3,-1}(e^n_{k,l}\otimes e^n_{3n+2-k,n-l})\\
&\quad=-(n+\tfrac53)e^{n+1}_{k+3,l+1}\otimes e^n_{3n+2-k,n-l}-L_{3,1}(n,k,l)e^{n-1}_{k,l}\otimes e^n_{3n+2-k,n-l}\\
&\qquad-(n+\tfrac53)e^n_{k,l}\otimes e^{n+1}_{3n+5-k,n+1-l}-L_{3,1}(n,3n+2-k,n-l)e^n_{k,l}\otimes e^{n-1}_{3n+2-k,n-l}.
\end{aligned}
\end{equation}
Terms with a missing level or weight are omitted. Multiplication by $(-1)^{k+l}\binom{3n+2}{k}\binom nl$ and summation gives the first-order action on $\Psi_n$.

Next apply $\mathcal Y_{3,1}$ to every term of this entire first-order expression. For example, its action on the first term of \eqref{eq:first-tensor-example} is
\begin{equation}\label{eq:second-tensor-example}
\begin{aligned}
&\mathcal Y_{3,1}\bigl(-(n+\tfrac53)e^{n+1}_{k+3,l+1}\otimes e^n_{3n+2-k,n-l}\bigr)\\
&\quad=(n+\tfrac53)(n+\tfrac83)e^{n+2}_{k+3,l+1}\otimes e^n_{3n+2-k,n-l}\\
&\qquad+(n+\tfrac53)L_{0,0}(n+1,k+3,l+1)e^n_{k,l}\otimes e^n_{3n+2-k,n-l}\\
&\qquad+(n+\tfrac53)^2e^{n+1}_{k+3,l+1}\otimes e^{n+1}_{3n+2-k,n-l}\\
&\qquad+(n+\tfrac53)L_{0,0}(n,3n+2-k,n-l)
 e^{n+1}_{k+3,l+1}\otimes e^{n-1}_{3n-1-k,n-1-l}.
\end{aligned}
\end{equation}
The second term, for instance, returns to the original tensor and has coefficient
\[
(n+\tfrac53)L_{0,0}(n+1,k+3,l+1)
=\frac{(k+3)^{\underline3}(l+1)}{27(n+2)(n+1)}.
\]
Applying the second operator to the other three terms completes this ordered product. Repeating these successive applications for the other three products in \eqref{eq:Q} gives their complete second-order expressions. Equations~\eqref{eq:path-table} and~\eqref{eq:cross-paths} record all resulting terms without suppressing any contribution.

Put $r'=3-r$, $t'=1-t$, and apply $Y'=Y_{3-2r',1-2t'}$ first and $Y=Y_{3-2r,1-2t}$ second. For the same factor the four paths in $\mathrm d\pi(Y)\mathrm d\pi(Y')e^n_{k,l}$ have coefficients
\begin{equation}\label{eq:path-table}
\begin{array}{c|c|c}
\text{first, second step}&\text{target}&\text{coefficient}\\\hline
\text{upper, upper}&e^{n+2}_{k+3,l+1}&U_nU_{n+1}\\
\text{upper, lower}&e^n_{k,l}&U_nL_{r,t}(n+1,k+r',l+t')\\
\text{lower, upper}&e^n_{k,l}&L_{r',t'}(n,k,l)U_{n-1}\\
\text{lower, lower}&e^{n-2}_{k-3,l-1}&L_{r',t'}(n,k,l)L_{r,t}(n-1,k-r,l-t).
\end{array}
\end{equation}
The last row is used only for $n\geq2$; the third only for $n\geq1$. The last coefficient simplifies to
\begin{equation}\label{eq:DD}
\frac{(k)^{\underline3}(3n+2-k)^{\underline3}l(n-l)}{D_nD_{n-1}},\quad n\geq2,
\end{equation}
independently of $(r,t)$.

Here is an explicit rule for every cross-factor term as well. For an existing input level let
\[
\begin{aligned}
c^+_{r,t}(n,k,l)&=-U_n,& s^+_{r,t}&=(r,t),\\
c^-_{r,t}(n,k,l)&=-L_{r,t}(n,k,l),& s^-_{r,t}&=(r-3,t-1).
\end{aligned}
\]
The minus option exists only for $n\geq1$. If $\mathbf k=(k,l)$ and $\mathbf j=(3n+2-k,n-l)$, the cross-factor part on $e^n_{\mathbf k}\otimes e^n_{\mathbf j}$ is exactly
\begin{align}\label{eq:cross-paths}
\sum_{\varepsilon,\eta\in\{+,-\}}\big[&
c^\varepsilon_{r,t}(n,\mathbf k)c^\eta_{r',t'}(n,\mathbf j)
e^{n+\varepsilon}_{\mathbf k+s^\varepsilon_{r,t}}\otimes e^{n+\eta}_{\mathbf j+s^\eta_{r',t'}}\nonumber\\
&+c^\varepsilon_{r',t'}(n,\mathbf k)c^\eta_{r,t}(n,\mathbf j)
e^{n+\varepsilon}_{\mathbf k+s^\varepsilon_{r',t'}}\otimes e^{n+\eta}_{\mathbf j+s^\eta_{r,t}}\big],
\end{align}
where $n+\varepsilon$ means $n+1$ or $n-1$ and absent options are omitted. Together with the two same-factor copies of \eqref{eq:path-table}, this gives the full sixteen terms of each ordered block. Multiply by $(-1)^{k+l}\binom{3n+2}k\binom nl$ and sum over the ranges in \eqref{eq:Psi}; then use the four signed products in \eqref{eq:Q}. This records the complete expansion without introducing singular low-level paths.

\subsection{Nonmatching sectors and the mixed contraction}\label{app:nonmatching}
The only potential sectors are $(n\pm2,n)$, $(n,n\pm2)$, $(n+1,n-1)$, $(n-1,n+1)$, $(n+1,n+1)$, $(n-1,n-1)$, and $(n,n)$. The first pair has common coefficient $U_nU_{n+1}$ multiplied by $-1+1+3-3=0$.

For the mixed sector $(n+1,n-1)$ the eight cross-factor paths, after the shifts of $k,l$ in \eqref{eq:cross-paths}, have common factor
\begin{equation}\label{eq:mixed-factor}
F_{n,k,l}=(-1)^{k+l}\frac{(3n+2)!}{(k-3)!(3n+2-k)!}\frac{n!}{(l-1)!(n-l)!}\frac{U_n}{D_n},
\end{equation}
for the valid target range $3\leq k\leq (3n+2)$, $1\leq l\leq n$, $n\geq1$. Their residual multipliers in the order of the eight terms in the aggregation below are
\begin{equation}\label{eq:mixed-cancel}
-1,-1,+1,+1,+3,+3,-3,-3,
\qquad\text{whose sum is zero}.
\end{equation}
The other mixed sector follows by interchanging tensor factors. Outside the valid range the target is absent; factorials outside that range are not evaluated. Compact invariance also rules out nonmatching sectors, but \eqref{eq:mixed-cancel} verifies the contraction explicitly.

\label{app:mixed}

All sums below give coefficients for $Q=\tfrac12\mathrm d\Pi(\Cas)$ on invariant vectors. For the detailed sums, use the weight notation
\[
E_{m,i}^{a,b}=\varphi_m(\xi_a^{3m+2}\boxtimes\xi_b^m),\quad i=1,2.
\]
The index $i$ identifies the tensor factor containing this matrix coefficient.

We use $\binom qj=0$ for nonnegative integers $q$ and $j\notin\{0,\ldots,q\}$. The retained target ranges are $0\leq k\leq3n+5$, $0\leq l\leq n+1$ for the upper block; $n\geq1$, $0\leq k\leq3n-1$, $0\leq l\leq n-1$ for the lower block; and $0\leq k\leq3n+2$, $0\leq l\leq n$ for the diagonal block. This is the scalar binomial convention for absent weights and summands; absent paths are omitted before coefficient evaluation.

For the mixed target
\[
E_{n+1,1}^{3n+5-2k,\,n+1-2l}\otimes E_{n-1,2}^{2k-3n-5,\,2l-n-1},
\]
the total coefficient is the following sum:
{\allowdisplaybreaks
\begin{align*}
&-(-1)^{k-3+l-1}\binom{3n+2}{k-3}\binom{n}{l-1}\rbra{n+\frac{5}{3}}\fs{(3n+5-k)^{\underline{3}}(n+1-l)}{9(3n+2)(n+1)^{\underline{2}}}\\
&-(-1)^{k+l}\binom{3n+2}{k}\binom{n}{l}\rbra{n+\frac{5}{3}}\fs{(k)^{\underline{3}}\,l}{9(3n+2)(n+1)^{\underline{2}}}\\
&+(-1)^{1+k-3+l}\binom{3n+2}{k-3}\binom{n}{l}\rbra{n+\frac{5}{3}}\fs{(3n+5-k)^{\underline{3}}\,l}{9(3n+2)(n+1)^{\underline{2}}}\\
&+(-1)^{1+k+l-1}\binom{3n+2}{k}\binom{n}{l-1}\rbra{n+\frac{5}{3}}\fs{(k)^{\underline{3}}(n+1-l)}{9(3n+2)(n+1)^{\underline{2}}}\\
&+3(-1)^{1+k-2+l-1}\binom{3n+2}{k-2}\binom{n}{l-1}\rbra{n+\frac{5}{3}}\fs{(3n+4-k)^{\underline{2}}(k-2)(n+1-l)}{9(3n+2)(n+1)^{\underline{2}}}\\
&+3(-1)^{1+k-1+l}\binom{3n+2}{k-1}\binom{n}{l}\rbra{n+\frac{5}{3}}\fs{(k-1)^{\underline{2}}(3n+3-k)l}{9(3n+2)(n+1)^{\underline{2}}}\\
&-3(-1)^{k-2+l}\binom{3n+2}{k-2}\binom{n}{l}\rbra{n+\frac{5}{3}}\fs{(3n+4-k)^{\underline{2}}(k-2)l}{9(3n+2)(n+1)^{\underline{2}}}\\
&-3(-1)^{k-1+l-1}\binom{3n+2}{k-1}\binom{n}{l-1}\rbra{n+\frac{5}{3}}\fs{(k-1)^{\underline{2}}(3n+3-k)(n+1-l)}{9(3n+2)(n+1)^{\underline{2}}}
\end{align*}
}
The binomial coefficients from \eqref{eq:Psi} combine with the downward falling factorials. For example,
\begin{align*}
&\binom{3n+2}{k-3}(3n+5-k)^{\underline3}
=\frac{(3n+2)!}{(k-3)!(3n+2-k)!},\\
&\binom n{l-1}(n+1-l)=\frac{n!}{(l-1)!(n-l)!}.
\end{align*}
The other seven shifted terms have the same factorial factors. By extracting the common multiplicative factor 
\[
(-1)^{k+l}\fs{(3n+2)!}{(k-3)!(3n+2-k)!}\fs{n!}{(l-1)!(n-l)!}\fs{n+\frac{5}{3}}{9(3n+2)(n+1)^{\underline{2}}},
\]
the eight residual multipliers are those in \eqref{eq:mixed-cancel}, whose sum is zero.
A similar argument applies to the component $E_{n-1,1}^{3n+5-2k,\,n+1-2l}\otimes E_{n+1,2}^{2k-3n-5,\,2l-n-1}$, which also vanishes.

For the downward shifts $(n-2,n)$ and $(n,n-2)$, the coefficient \eqref{eq:DD} is independent of $(r,t)$ and again has signed multiplier $-1+1+3-3=0$. This calculation is used only for $n\geq2$; such paths are absent at the two lowest levels.

\subsection{Aggregation of the retained coefficients}\label{app:aggregation}

The coefficient for\\
$E_{n+1,1}^{3n+5-2k,\,n+1-2l}\otimes E_{n+1,2}^{2k-3n-5,\,2l-n-1}$ is
\begin{align*}
\rbra{n+\frac{5}{3}}^2&\Bigg[ (-1)^{1+k-3+l-1}\binom{3n+2}{k-3}\binom{n}{l-1}+(-1)^{1+k+l}\binom{3n+2}{k}\binom{n}{l}\\
&\quad+(-1)^{k-3+l}\binom{3n+2}{k-3}\binom{n}{l}+(-1)^{k+l-1}\binom{3n+2}{k}\binom{n}{l-1}\\
&\quad+3(-1)^{k-2+l-1}\binom{3n+2}{k-2}\binom{n}{l-1}+3(-1)^{k-1+l}\binom{3n+2}{k-1}\binom{n}{l}\\
&\quad -3(-1)^{k-2+l}\binom{3n+2}{k-2}\binom{n}{l}-3(-1)^{k-1+l-1}\binom{3n+2}{k-1}\binom{n}{l-1}\Bigg].
\end{align*}
Factoring out the term $-(-1)^{k+l}\rbra{n+\fs{5}{3}}^2\fs{(3n+2)!}{k!(3n+5-k)!}\fs{n!}{l!(n-l+1)!}$, the remaining expression simplifies to:
\begin{align*}
&(k)^{\underline{3}}\,l+(3n+5-k)^{\underline{3}}(n-l+1)+(k)^{\underline{3}}(n+1-l)+(3n+5-k)^{\underline{3}}\,l\\
&+3(k)^{\underline{2}}(3n+5-k)l+3k(3n+5-k)^{\underline{2}}(n+1-l)\\
&+3(k)^{\underline{2}}(3n+5-k)(n+1-l)+3k(3n+5-k)^{\underline{2}}\,l\\
&=\rbra{\sum_{m=0}^3\binom{3}{m}(k)^{\underline{m}}(3n+5-k)^{\underline{3-m}}}(n+1) = (3n+5)^{\underline{3}}(n+1).
\end{align*}
The falling-factorial Vandermonde identity, for positive integers $a, b$ and $n$, is
\[
\sum_{k=0}^n\binom{n}{k} (a)^{\underline{k}}(b)^{\underline{n-k}}=(a+b)^{\underline{n}}.
\]
The summed coefficient is therefore
\[
-(-1)^{k+l}\binom{3n+5}{k}\binom{n+1}{l}\times \rbra{n+\frac{5}{3}}^2.
\]
Thus the coefficient of $\Psi_{n+1}$ in the half-Casimir $Q$ is $-\rbra{n+\frac{5}{3}}^2$.

The same upper coefficient can be collected directly in binomial form. Each of the eight contributing terms raises both levels. After these index shifts, their coefficient is
\begin{equation}\label{eq:upper-aggregation-main}
-(-1)^{k+l}(n+\tfrac53)^2
\sum_{j=0}^3\binom3j\binom{3n+2}{k-j}
\left(\binom n{l-1}+\binom nl\right).
\end{equation}
The weights $1,3,3,1$ in this four-term sum come from the four signed products after combining the two possible second-factor shifts; it is precisely the eight-term upper block, grouped in pairs. With out-of-range binomial coefficients set to zero, Pascal's identity and binomial convolution give
\[
\binom n{l-1}+\binom nl=\binom{n+1}l,
\quad
\sum_{j=0}^3\binom3j\binom{3n+2}{k-j}=\binom{3n+5}k.
\]
Thus \eqref{eq:upper-aggregation-main} is $-(n+\frac53)^2$ times the coefficient of this tensor in $\Psi_{n+1}$. This proves the upper coefficient of $Q\Psi_n$.

\paragraph{Lower matching coefficient.}
For $n\geq1$, the lower target is
\[
E_{n-1,1}^{3n-1-2k,\, n-1-2l}\otimes E_{n-1,2}^{2k+1-3n,\, 2l+1-n}.
\]
The eight terms, after their index shifts, are
{\allowdisplaybreaks
\begin{align*}
    &-(-1)^{k+l}\binom{3n+2}{k}\binom{n}{l}\rbra{\fs{(3n+2-k)^{\underline{3}}(n-l)}{9(3n+2)(n+1)^{\underline{2}}}}^2 \\
    &-(-1)^{k+3+l+1}\binom{3n+2}{k+3}\binom{n}{l+1}\rbra{\fs{(k+3)^{\underline{3}}(l+1)}{9(3n+2)(n+1)^{\underline{2}}}}^2 \\
    &+(-1)^{k+l+1}\binom{3n+2}{k}\binom{n}{l+1}\rbra{\fs{(3n+2-k)^{\underline{3}}(l+1)}{9(3n+2)(n+1)^{\underline{2}}}}^2 \\
    &+(-1)^{k+3+l}\binom{3n+2}{k+3}\binom{n}{l}\rbra{\fs{(k+3)^{\underline{3}}(n-l)}{9(3n+2)(n+1)^{\underline{2}}}}^2 \\
    &+3(-1)^{k+1+l}\binom{3n+2}{k+1}\binom{n}{l}\rbra{\fs{(3n+1-k)^{\underline{2}}(k+1)(n-l)}{9(3n+2)(n+1)^{\underline{2}}}}^2 \\
    &+3(-1)^{k+2+l+1}\binom{3n+2}{k+2}\binom{n}{l+1}\rbra{\fs{(k+2)^{\underline{2}}(3n-k)(l+1)}{9(3n+2)(n+1)^{\underline{2}}}}^2 \\
    &-3(-1)^{k+1+l+1}\binom{3n+2}{k+1}\binom{n}{l+1}\rbra{\fs{(3n+1-k)^{\underline{2}}(k+1)(l+1)}{9(3n+2)(n+1)^{\underline{2}}}}^2 \\
    &-3(-1)^{k+2+l}\binom{3n+2}{k+2}\binom{n}{l}\rbra{\fs{(k+2)^{\underline{2}}(3n-k)(n-l)}{9(3n+2)(n+1)^{\underline{2}}}}^2.
\end{align*} 
}
Extracting the common factor 
\[
    \mathcal{C}_{n,k,l} \coloneqq (-1)^{k+l+1}\fs{(3n+2)!}{k!(3n-1-k)!}\fs{n!}{l!(n-1-l)!}\rbra{\fs{1}{9(3n+2)(n+1)^{\underline{2}}}}^2,
\]
the remaining combinatorial sum simplifies as follows:
\begin{align*}
    &(3n+2-k)^{\underline{3}}(n-l) + (k+3)^{\underline{3}}(l+1) + (3n+2-k)^{\underline{3}}(l+1) + (k+3)^{\underline{3}}(n-l) \\
    &+\;
    3(3n+1-k)^{\underline{2}}(k+1)(n-l) + 3(k+2)^{\underline{2}}(3n-k)(l+1) \\
    &+\; 3(3n+1-k)^{\underline{2}}(k+1)(l+1) + 3(k+2)^{\underline{2}}(3n-k)(n-l) \\
    &=\; \rbra{\sum_{m=0}^3\binom{3}{m}(k+1)_m(3n-k)_{3-m}}(n-l+l+1) = (3n+1)_3(n+1).
\end{align*}
The rising-factorial Vandermonde identity, for positive integers $a, b$ and $n$, is
\[
    \sum_{j=0}^n\binom{n}{j} (a)_j(b)_{n-j}=(a+b)_{n}.
\]
The summed coefficient is therefore
\begin{align*}
    \mathcal{C}_{n,k,l} \cdot (3n+1)_3(n+1) 
    &=
    -(-1)^{k+l}\binom{3n-1}{k}\binom{n-1}{l}\fs{(3n)_3(3n+1)_3\,n(n+1)}{(3n+2)^2(3n+3)^2(3n)^2} \\
    &= -\rbra{n+\frac{1}{3}}^2 (-1)^{k+l}\binom{3n-1}{k}\binom{n-1}{l}.
\end{align*}
Thus the coefficient of $\Psi_{n-1}$ in the half-Casimir $Q$ is $-\rbra{n+\frac{1}{3}}^2$.

Equivalently, the lower coefficient before division by the invariant tensor coefficient is
\begin{equation}\label{eq:lower-aggregation-main}
\begin{aligned}
&-(-1)^{k+l}\frac{(3n+2)!}{k!(3n-1-k)!}
\frac{n!}{l!(n-1-l)!}\frac{1}{[9(3n+2)n(n+1)]^2}\\
&\quad\times(n+1)\sum_{j=0}^3\binom3j(k+1)_j(3n-k)_{3-j}.
\end{aligned}
\end{equation}
Here $0\leq k\leq3n-1$ and $0\leq l\leq n-1$. The factor $n+1$ is the sum $(n-l)+(l+1)$ from the two second-factor shifts. The remaining four terms are the rising-factorial convolution
\[
\sum_{j=0}^3\binom3j(k+1)_j(3n-k)_{3-j}=(3n+1)_3.
\]
Dividing \eqref{eq:lower-aggregation-main} by the coefficient $(-1)^{k+l}\binom{3n-1}k\binom{n-1}l$ in $\Psi_{n-1}$ leaves
\[
-\frac{(3n)_3(3n+1)_3\,n(n+1)}{[9(3n+2)n(n+1)]^2}
=-\frac{(3n+1)^2}{9}
=-(n+\tfrac13)^2.
\]
This gives the lower coefficient of $Q\Psi_n$.

\paragraph{Diagonal matching coefficient.}
The diagonal tensor component is
\[
E_{n,1}^{3n+2-2k,n-2l}\otimes E_{n,2}^{2k-3n-2,2l-n}.
\] 
The upper--lower and lower--upper paths give two groups, with denominators involving $(n+2)^{\underline{2}}$ and $(n+1)^{\underline{2}}$, respectively.

The first group has denominator $27(n+2)^{\underline{2}}$:
{\allowdisplaybreaks
\begin{align*}
    &(-1)^{k+l+1}\binom{3n+2}{k}\binom{n}{l}\Bigg[
    \fs{(k+3)^{\underline{3}}(l+1)}{27(n+2)^{\underline{2}}}+\fs{(3n+5-k)^{\underline{3}}(n+1-l)}{27(n+2)^{\underline{2}}} \\
    &\quad+\fs{(k+3)^{\underline{3}}(n+1-l)}{27(n+2)^{\underline{2}}}+\fs{(3n+5-k)^{\underline{3}}(l+1)}{27(n+2)^{\underline{2}}} \\
    &\quad+3\fs{(3n+3-k)(k+2)^{\underline{2}}(l+1)}{27(n+2)^{\underline{2}}}+3\fs{(3n+4-k)^{\underline{2}}(k+1)(n+1-l)}{27(n+2)^{\underline{2}}} \\
    &\quad+3\fs{(3n+3-k)(k+2)^{\underline{2}}(n+1-l)}{27(n+2)^{\underline{2}}}+3\fs{(3n+4-k)^{\underline{2}}(k+1)(l+1)}{27(n+2)^{\underline{2}}}\Bigg].
\end{align*}
}

The rising-factorial Vandermonde identity gives:
\begin{align*}
    &(-1)^{k+l+1}\binom{3n+2}{k}\binom{n}{l}\rbra{\sum_{m=0}^3\binom{3}{m}(3n+3-k)_m(k+1)_{3-m}}\!\fs{(n+1-l+l+1)}{27(n+2)^{\underline{2}}} \\
    &\quad= -(-1)^{k+l}\binom{3n+2}{k}\binom{n}{l}(3n+4)_3\fs{(n+2)}{27(n+2)(n+1)} \\
    &\quad= (-1)^{k+l}\binom{3n+2}{k}\binom{n}{l}\cdot\rbra{-\fs{\rbra{n+\frac{4}{3}}\rbra{n+\frac{5}{3}}(n+2)}{n+1}}.
\end{align*}

For $n\geq1$ only, the second component consists of terms with the denominator $27(n+1)^{\underline{2}}$:
{\allowdisplaybreaks
\begin{align*}
    (-1)^{k+l+1}\binom{3n+2}{k}\binom{n}{l}\Bigg[&
    \fs{(3n+2-k)^{\underline{3}}(n-l)}{27(n+1)^{\underline{2}}}+\fs{(k)^{\underline{3}}l}{27(n+1)^{\underline{2}}} \\
    &\quad+\fs{(3n+2-k)^{\underline{3}}l}{27(n+1)^{\underline{2}}}+\fs{(k)^{\underline{3}}(n-l)}{27(n+1)^{\underline{2}}} \\
    &\quad+3\fs{(3n+2-k)^{\underline{2}}k(n-l)}{27(n+1)^{\underline{2}}}+3\fs{(k)^{\underline{2}}(3n+2-k)l}{27(n+1)^{\underline{2}}} \\
    &\quad+3\fs{(3n+2-k)^{\underline{2}}kl}{27(n+1)^{\underline{2}}}+3\fs{(k)^{\underline{2}}(3n+2-k)(n-l)}{27(n+1)^{\underline{2}}}\Bigg].
\end{align*}
}
As in the $(n+1, n+1)$ case, the falling-factorial Vandermonde identity reduces this expression to
\begin{align*}
    & (-1)^{k+l+1}\binom{3n+2}{k}\binom{n}{l}\rbra{\sum_{m=0}^3\binom{3}{m}(3n+2-k)^{\underline{m}}(k)^{\underline{3-m}}}\fs{(n-l+l)}{27(n+1)^{\underline{2}}} \\
    &\quad= -(-1)^{k+l}\binom{3n+2}{k}\binom{n}{l}(3n+2)^{\underline{3}}\fs{n}{27(n+1)n} 
\end{align*}
\[
= (-1)^{k+l}\binom{3n+2}{k}\binom{n}{l}\cdot \rbra{-\fs{n\rbra{n+\frac{1}{3}}\rbra{n+\frac{2}{3}}}{n+1}}.
\]

Combining the two groups gives the coefficient of $\Psi_{n}$:
\[
    -\fs{\rbra{n+\frac{4}{3}}\rbra{n+\frac{5}{3}}(n+2)}{n+1}-\fs{n\rbra{n+\frac{1}{3}}\rbra{n+\frac{2}{3}}}{n+1}
\]
to the half-Casimir $Q$.
For $n=0$ the second diagonal group is absent, so only the regular upper--lower contribution is used. The three aggregated coefficients belong to $Q$. Doubling them, as required by \eqref{eq:Q}, proves \eqref{equ:4.recurrence_final}.

\subsection{Separate boundary calculations}\label{app:boundary}
At level zero use only the upper term of \eqref{eq:one-factor}; a subsequent lower action from level one is regular. At level one, a downward step reaches zero, so a second downward step in that factor is absent. A downward step in the other level-one factor remains allowed. No $D_0$ or $D_{-1}$ is evaluated.

The following table records the four signed ordered contributions to the tensor component $k=l=0$ of each retained output level. That coefficient in $\Psi_j$ is one. Together with the invariant-line construction and the nonmatching cancellations, it determines the whole vector. Columns $Q_1,Q_2,Q_3,Q_4$ refer in order to the four signed products in \eqref{eq:Q}.
\[
\renewcommand{\arraystretch}{1.6}
\begin{array}{cc|rrrr|r}
 n&\text{output}&Q_1&Q_2&Q_3&Q_4&Q\\\hline
0&1&-\frac{25}{9}&0&0&0&-\frac{25}{9}\\
0&0&-\frac{11}{9}&-\frac{11}{9}&-1&-1&-\frac{40}{9}\\
1&2&-\frac{64}{9}&0&0&0&-\frac{64}{9}\\
1&1&-\frac{143}{27}&-\frac{58}{27}&-\frac{16}{9}&-\frac{11}{9}&-\frac{94}{9}\\
1&0&-\frac{22}{45}&-\frac{22}{45}&-\frac{2}{5}&-\frac{2}{5}&-\frac{16}{9}.
\end{array}
\]
These regular one-step sums give exactly \eqref{eq:casimir-boundary} after multiplication by two.

\section*{Acknowledgements}

This article is based on material from the author's doctoral dissertation. OpenAI's ChatGPT was used as an organizational aid in extracting the relevant material from the dissertation, restructuring it into article form, and assisting with linguistic and presentational revisions. The mathematical content, including all statements, proofs, and their final verification, remains the responsibility of the author.

\bibliographystyle{amsplain}
\bibliography{references}
\end{document}